\documentclass[english]{amsart}

\usepackage{amsmath,amssymb,amsthm,array,multicol,verbatim,graphicx}
\usepackage{epsfig,float}
\usepackage{psfrag}
\usepackage{url}
\usepackage[text={445pt,630pt}, centering, a4paper,
marginparwidth=2cm]{geometry}
\usepackage{tikz}
\usepackage[shortlabels]{enumitem}
\usepackage[english]{babel}
\usepackage[utf8]{inputenc}
\usepackage[T1]{fontenc}		
\usepackage{tabularx}
\usepackage[hidelinks,pdfusetitle]{hyperref}

\theoremstyle{plain}
\newtheorem{thrm}{Theorem}[section]
\newtheorem{crllr}[thrm]{Corollary}
\newtheorem{prpstn}[thrm]{Proposition}
\newtheorem{lmm}[thrm]{Lemma}
\newtheorem{dfntn}[thrm]{Definition}

\theoremstyle{definition}
\newtheorem{rmrk}[thrm]{Remark}

\newcommand{\reff}[1]{(\ref{#1})}

\newcommand{\ca}{{\mathcal A}}

\newcommand{\ce}{{\mathcal E}}

\newcommand{\cf}{{\mathcal F}}

\newcommand{\ch}{{\mathcal H}}

\newcommand{\cl}{{\mathcal L}}

\newcommand{\cp}{{\mathcal P}}

\newcommand{\cs}{{\mathcal S}}

\newcommand{\cu}{{\mathcal U}}

\newcommand{\cy}{{\mathcal Y}}
\newcommand{\cz}{{\mathcal Z}}

\newcommand{\A}{{\mathbb A}}

\newcommand{\E}{{\mathbb E}}

\newcommand{\N}{{\mathbb N}}
\renewcommand{\P}{{\mathbb P}}

\newcommand{\T}{{\mathbb T}}
\newcommand{\Z}{{\mathbb Z}}

\newcommand{\rA}{{\rm A}}

\newcommand{\rP}{{\rm P}}
\newcommand{\rE}{{\rm E}}

\newcommand{\bt}{{\mathbf t}}
\newcommand{\bs}{{\mathbf s}}
\newcommand{\bv}{{\mathbf v}}

\newcommand{\ind}{{\bf 1}}

\renewcommand{\root}{{\varnothing}}

\newcommand{\supp}{{\rm supp}\;}

\newcommand{\Card}{{\rm Card}\;}

\newcommand{\val}[1]{\mathop{\left| #1 \right|}\nolimits}
\newcommand{\inv}[1]{\mathop{\frac{1}{ #1}}\nolimits}
\newcommand{\expp}[1]{\mathop {\mathrm{e}^{ #1}}}

\newcommand{\mt}{{\breve \tau}}

\usepackage[draft]{fixme} 
\FXRegisterAuthor{jf}{ajf}{JFD}
\FXRegisterAuthor{ra}{ara}{RA}
\fxusetheme{colorsig}

\begin{document}

\title{An introduction to Bienaym\'e-Galton-Watson trees and their local limits}
\date{\today}
 
\author{Romain Abraham}
\address{
Institut Denis Poisson
Universit\'e d'Orl\'eans, Universit\'e de Tours, CNRS,
B.P. 6759,
45067 Orl\'eans cedex 2,
France.}
\author{Jean-Fran{\c c}ois Delmas}
\address{CERMICS, \'Ecole des Ponts, France.}



\begin{abstract} 
The aim of this lecture is to give an overview of old and new results
on Bienaym\'e-Galton-Watson (BGW) trees. After introducing the framework of discrete
trees, we first give alternative proofs of classical results on the
extinction probability of BGW processes and on the
description of the processes conditioned on extinction or on
non-extinction. Then, we study recent  local limits of critical or
sub-critical BGW trees conditioned to be large.
\end{abstract}

%
%

\maketitle
\tableofcontents

\section{Introduction}

The first draft  of those notes has  been written for a  course given at
Hamamet in 2014, it has then been completed with more examples and up to
date references. It concerns the so-called Bienaym\'e-Galton-Watson (BGW
for short) process which can be considered as the first stochastic model
for population evolution.   It was named after  the French mathematician
I.-J.~Bienaym\'e  (1845)  and  the   British  scientists  F.~Galton  and
H.~W.~Watson  (1874) who  studied it.   I.-J.~Bienaym\'e considered  the
probability  of  the  male  population extinction  in  \cite{b:gw};  its
communication, which  is reproduced in~\cite{k:history},  indicates that
he  knew  the  right  answer (see  also  the  study~\cite{bru}  on  the
Bienaym\'e's proof  of the criticality theorem  for branching processes).
Later on and independently, F.~Galton, who was studying human evolution, published in 1873
in Educational Times a question on  the probability of extinction of the
noble surnames in the UK. It was a very short communication which can be
copied integrally here:
\begin{quote}
``PROBLEM 4001: A large nation, of whom we will only concern ourselves with adult males, $N$ in number, and who each bear separate surnames colonise a district. Their law of population is such that, in each generation, $a_0$ per cent of the adult males have no male children who reach adult life; $a_1$ have one such male child; $a_2$ have two; and so on up to $a_5$ who have five.
Find (1) what proportion of their surnames will have become extinct after $r$ generations; and (2) how many instances there will be of the surname being held by $m$ persons.''
\end{quote}
In more  modern terms,  he supposes that  all the  individuals reproduce
independently  from  each  others  and   have  all  the  same  offspring
distribution. After  receiving no valuable  answer to that  question, he
directly contacted H.~W.~Watson  and worked together on  the problem. They
published an article one year later \cite{gw:pef} where they proved that
the  probability  of extinction  is  a  fixed  point of  the  generating
function  of the  offspring  distribution (which  is  true, see  Section
\ref{sec:ce}) and concluded  a bit too rapidly that  this probability is
always equal to  1 (which is false, see also  Section \ref{sec:ce}). For
further historical  comments on  BGW processes,  we refer  to D.~Kendall
\cite{k:history} for the ``Genealogy of genealogy branching process'' up
to           1975           as           well           as           the
Lecture\footnote{\url{http://www.math.chalmers.se/~jagers/Branching History.pdf}}
at the Oberwolfach  Symposium on ``Random Trees'' in  2009 by P.~Jagers.
In order to track the genealogy of  the population of a BGW process, one
can consider  the so called  genealogical trees  or BGW trees,  which is
currently  an  active  domain  of   research.   We  refer  to  T.~Harris
\cite{h:tbp} and  K.~Athreya and P.~Ney \cite{an:bp}  for most important
results on BGW processes, to  M.~Kimmel and D.~Axelrod \cite{ka:bpb} and
P.~Haccou, P.~Jagers  and V.~Vatutin  \cite{hjv:pb} for  applications in
biology,  to M.~Drmota  \cite{d:rt} and  S.~Evans \cite{e:rt}  on random
discrete trees including BGW trees (see also J.~Pitman \cite{p:csp} on a
more   combinatorial   aspect   and  T.~Duquesne   and   J.-F.~Le   Gall
\cite{dlg:rt}  for  scaling  limits  of  BGW trees  which  will  not  be
presented here).

\medskip We introduce in the first  chapter  the framework
of   discrete  random   trees,  which   may  be   attributed  to   Neveu
\cite{n:apghw}.  We then  use this framework to construct  BGW trees that
describe the genealogy of  a BGW process. It is very  easy to recover the
BGW process from the  BGW tree as it is just the  number of individuals at
each generation. We then give alternative proofs of classical results on
BGW processes  using the tree  formalism. We  focus in particular  on the
extinction probability (which  was the first question  of F.~Galton) and
on  the description  of the  processes conditioned  on extinction  or
non extinction. 

\medskip In  a second chapter, we  focus on local limits  of conditioned
BGW trees. In  the critical and sub-critical cases (these  terms will be
explained in  the first chapter),  the population becomes  a.s.\ extinct
and the associated genealogical tree is  finite. However, it has a small
but  positive probability  of  being  large (this  notion  must be  made
precise). The  question that arises is  to describe the law  of the tree
conditioned  of being  large,  and  to say  what  exceptional event  has
occurred  so that  the  tree is  not  typical. A  first  answer to  this
question is due to H.~Kesten~\cite{k:sbrwrc} who conditioned a BGW tree
to reach height $n$ and look at the limit in distribution when $n$ tends
to infinity. There  are however other ways of conditioning  a tree to be
large:      conditioning     on      having     many      nodes,     see
S.~Janson~\cite{j:sgtcgwrac}  and  references  therein; on  having  many
leaves, see I. Kortchemski~\cite{kor:GW-cond-leaf}; on having  nodes with large degrees,  see X. He~\cite{h:cgwtmod}
and  B.~Stufler~\cite{stufler20};   on  having  large   population,  see
R.~Abraham and J.-F.~Delmas~\cite{ad:apegw};...  In most of those cases,
the local limit is an infinite random tree with either an infinite spine
(the so  called Kesten's tree);  or with  an infinite backbone,  see for
example~\cite{abd:nlcgwt,ad:apegw}; or  with a  node of  infinite degree
(the so called condensation phenomenon),  see for example T.~Jonsson and
S.~O.~Stef\'ansson~\cite{js:cnt},     S.~Janson~\cite{j:sgtcgwrac}    or
B.~Stufler~\cite{stufler19}  and the  references therein.   The approach
presented    here   is    mainly    based   on    our   previous    work
\cite{ad:llcgwtisc,ad:llcgwtcc}.

\section{Bienaym\'e-Galton-Watson  trees and extinction} 

We  intend to  give a  short introduction  to Bienaym\'e-Galton-Watson  (BGW) trees,
which  is  an  elementary  model  for  the  genealogy  of  a  branching
population. The  BGW process, which can  be defined directly from  the BGW
tree,   describes   the  evolution   of   the   size  of   a   branching
population.  Roughly speaking,  each  individual of  a given  generation
gives birth to  a random number of children in  the next generation. The
distribution probability  of the random  number of children,  called the
offspring  distribution,  is the  same  for  all the  individuals.   The
offspring   distribution    is   called   sub-critical,    critical   or
super-critical if its  mean is respectively strictly less  than 1, equal
to 1, or strictly greater than 1.

We denote by $\N=\Z_+$ the set of non-negative integers and by
$\N^*=\N\setminus\{0\}$ the set of positive integers.

\medskip

We describe more precisely the BGW process. Let  $\zeta$ be  a random variable
taking   values  in   $\N$  with   distribution  $p=(p(n),   n\in  \N)$:
$p(n)=\P(\zeta=n)$.   We denote  by $m=\E[\zeta]$  the mean  of $\zeta$.
Let  $g(r)=\sum_{k\in  \N}  p(k)\,  r^k=\E\left[r^\zeta\right]$  be  the
generating  function of  $p$ defined  on  $[0,1]$.  We  recall that  the
function $g$ is convex, with $g'(1)=\E[\zeta]\in [0, +\infty ]$.

The  BGW  process $Z=(Z_n,  n\in  \N)$  with offspring  distribution  $p$
describes the evolution of the size of a population issued from a single
individual under the following assumptions:
\begin{itemize}
   \item $Z_n$ is the size of the population at time or generation
     $n$. In particular, $Z_0=1$. 
   \item Each individual alive at time  $n$ dies at generation $n+1$ and
     gives birth to a random number of children at time $n+1$, which is
     distributed as $\zeta$ and independent of the number of children of
     other individuals. 
\end{itemize}
 We can
define the process $Z$ more formally. Let $(\zeta_{i, n}; i\in \N, n\in
\N)$ be independent random variables  distributed as $\zeta$.
We set $Z_0=1$ and, with the convention $\sum_\emptyset=0$,  for $n\in \N^*$:
\begin{equation}
   \label{eq:def-Z}
Z_{n}=\sum_{i=1}^{Z_{n-1}} \zeta_{i, n}.
\end{equation}
Here  the  random  variable   $\zeta_{i,n}$  represents  the  number  of
offspring of the $i$-th individual alive at generation $n$.

The genealogical tree, or BGW tree, associated with the BGW process will be
described in Section \ref{sec:GW} after an introduction to discrete
trees given in Section \ref{sec:discrete}.

We say that
the population is extinct at time $n$ if $Z_n=0$ (notice that it is then
extinct at any further time, and thus $\{Z_n=0\}\subset \{Z_{n+1}=0\}$). 
The extinction event $\ce$ corresponds to:
\begin{equation}
   \label{eq:def-E}
\ce=\{\exists n\in \N \text{ s.t. } Z_n=0\}= \bigcup_{n\in \N}  \{Z_n=0\}.
\end{equation}

We  shall  compute  the  extinction  probability  $\P(\ce)$  in  Section
\ref{sec:ce}  using the  BGW  tree  setting (we  stress  that the  usual
computation  relies  on  the  properties of  $Z_n$  and  its  generating
function),    see    Corollary    \ref{cor:sub-critical}    and    Lemma
\ref{lem:sur-critique} which  state that $\P(\ce)$ is  the smallest root
of  $g(r)=r$ in  $[0,1]$. In  particular the  extinction is  almost sure
(a.s.) in the sub-critical case and critical case (unless $p(1)=1$). The
advantage  of the  proof provided  in Section  \ref{sec:ce}, is  that it
directly provides  the distribution of  the super-critical BGW  tree and
process   conditionally    on   the   extinction   event,    see   Lemma
\ref{lem:sur-critique}.

In  Section \ref{sec:super-cec},  we  describe the  distribution of  the
super-critical  BGW  tree  conditionally  on  the non-extinction  event,  see
Corollary   \ref{cor:TS}.  In Section   \ref{sec:asymptZ},  we   study 
asymptotics of the BGW process in the super-critical case, see
Theorem  \ref{theo:ks}.  We  prove  this result  from  Kesten  and  Stigum
\cite{ks:ltmgwp}  by following  the  proof  of Lyons,  Pemantle  and  Peres
\cite{lpp:cp}, which  relies on a change of  measure on  the genealogical
tree (this proof is also clearly exposed in Alsmeyer's lecture notes\footnote{\url{http://wwwmath.uni-muenster.de/statistik/lehre/WS1011/SpezielleStochastischeProzesse/}}). In particular we shall use Kesten's tree which is an elementary
multi-type BGW tree. It is defined in Section \ref{sec:defK} and it will
play a central role in Chapter \ref{sec:loc}.

\subsection{The set of discrete trees}
\label{sec:discrete}
We recall Neveu's formalism \cite{n:apghw} for ordered rooted trees. We set:
\[
\cu=\bigcup _{n\ge 0}{(\N^*)^n}
\]
the set  of finite  sequences of positive  integers with  the convention
$(\N^*)^0=\{\root\}$.  The  set $\cu$  is sometimes called  Ulam's tree.
For $n\geq  1$ and $u=(u_1,  \ldots, u_n) \in  \cu$, we set  $H(u)=n$ the
height or generation of $u$ and set $H(\root)= 0$. If $u$ and $v$ are two
sequences  of $\cu$,  we denote  by $uv$  the concatenation  of the  two
sequences, with the convention that $uv=vu=v$ if $u=\root$.  We define a
partial   order   on   $\cu$   called   the   genealogical   order   by:
$v\preccurlyeq u$  if there exists $w\in  \cu$ such that $u=vw$.  We say
that $v$ is an ancestor of $u$ and write $v\prec u$ if $v\preccurlyeq u$
and $v\neq u$.  The set of ancestors of $u\in \cu$ is the set:
\begin{equation}
   \label{eq:Au}
   \rA_u=\{v\in \cu; v\prec u\}.
\end{equation}
We set  $\bar \rA_u=\rA_u\cup\{u\}$  and notice  that $\root  \in \rA_u$
unless $u$  is the root.   The most recent  common ancestor of  a subset
$\bs$ of  $ \cu$, denoted  by $\text{MRCA}(\bs)$, is the  unique element
$v$ of $\bigcap_{u\in \bs} \bar\rA_u$  with maximal height. For the MRCA
of  two nodes,  say  $u$ and  $v$,  we  simply write  $v  \wedge u$  for
$\text{MRCA}(\{v,u\})$.  We  consider the lexicographic order  on $\cu$:
for  $u,v\in\cu$, we  set $v<u$  if either  $v\prec u$  or ($v=wjv'$  and
$u=wiu'$) with $w=u\wedge v$, and $j<i$ for some $i,j\in\N^*$.

\bigskip

A tree $\bt$ is a subset of $\cu$ that satisfies:
\begin{itemize}
\item $\root\in\bt$,
\item If  $u\in\bt$ and $v \prec u$, then $v\in  \bt$. 
\item For every $u\in \bt$, there exists 
  $k_u(\bt)\in \N\cup\{+\infty \}$ such that, for every  $i\in \N^*$,  $ui\in \bt$ iff $1\leq i\leq k_u(\bt)$. 
\end{itemize}

The integer $k_u(\bt)$ represents the number of offsprings of the node
$u\in \bt$.  The node $u\in \bt$  is called a leaf if $k_u(\bt)=0$ and
it is said infinite if $k_u(\bt)=+\infty $.  By convention, we shall set
$k_u(\bt)=-1$ if $u\not\in  \bt$.  The node $\root$  is called the
root of $\bt$.  A finite tree is represented in Fig.~\ref{fig:treeT0}.

\begin{figure}[H]
\begin{tikzpicture}
\node {$\root$} [grow=north]
child {node{(3)}
	child{node{(3,1)}
		child{node{(3,1,2)}}
		child{node{(3,1,1)}}
		}
	}
child {node {(2)}}
child {node {(1)}
	child{node{(1,3)}}
	child{node{(1,2)}}
	child{node{(1,1)}}
	}
;
\end{tikzpicture}
\caption{A finite tree $\bt$.}
\label{fig:treeT0}
\end{figure}


We   denote  by   $\T_\infty $  the   set  of   trees and by $\T$ the
subset of trees with no infinite node: 
\[
\T=\{\bt\in  \T_\infty ;\, k_u(\bt)<+\infty, \, \forall u\in \bt \}.
\]

Let us stress that the offspring of one individual are ordered;  this
amounts to  consider  planar trees. In particular the  two trees of
Fig.~\ref{fig:planar} are different.

\begin{figure}[H]
\begin{center}
\begin{tikzpicture}
\node {$\root$} [grow=north]
child {node {(2)}
	child{node{(2,1)}}
	}
child {node {(1)}
	}
;
\end{tikzpicture}
\hspace{2cm}
\begin{tikzpicture}
\node {$\root$} [grow=north]
child {node {(2)}
	}
child {node {(1)}
	child{node{(1,1)}}
	}
;
\end{tikzpicture}
\end{center}
\caption{Two different planar trees.}
\label{fig:planar}
\end{figure}

Let  $\bt\in \T_\infty $. We set $\sharp \bt=\Card(\bt)$ and notice that:
\begin{equation}\label{eq:sum_k}
\sum_{u\in\bt}k_u(\bt)=\sharp \bt -1.
\end{equation}
The set of its leaves is  $\cl_0(\bt)=\{u\in \bt;
k_u(\bt)=0\}$. Its height and its width at level $h\in \N$  are
respectively  defined by:
\[
H(\bt)=\sup\{H(u);\, u\in\bt\} 
\quad\text{and}\quad
z_h(\bt) =\Card \left(\{u\in \bt;\,  (u)=h\}\right);
\]
they  can be infinite.
Notice that  $\bt\in
\T$ if and only if  $z_h(\bt)$ is  finite for all $h\in \N$.  
For $u\in \bt$, we define the 
sub-tree  $\cs_u(\bt)$ of $\bt$ ``above'' $u$ as:
\begin{equation}
   \label{eq:defSu}
\cs_u(\bt)=\{v\in\cu;\ uv\in\bt\}=\{v\in\cu;\ u\in\bar A_v\}.
\end{equation}

We   will   mainly   consider   trees   in   $\T$,   but   for   Section
\ref{sec:non-gene} where we shall consider trees with one infinite node.
We  denote  by $\T_0$  the  countable  subset  of  finite trees  and  by
$\T^{(h)}\subset \T_0$ the subset of  finite trees with height less than
$h\in \N$:
\begin{equation}
   \label{eq:t0-th}
   \T_0=\{\bt\in  \T ;\,\sharp \bt<+\infty \}
   \quad\text{and}\quad
\T^{(h)}=\{\bt \in \T; \, H(\bt)\leq h\}.
\end{equation}
For $v=(v_k, k\in \N^*)\in  (\N^*)^{\N^*} $, we set $\bar v_n=(v_1,
\ldots,   v_n)$  for  $n\in   \N$,  with   the  convention   that  $\bar
v_0=\root$ and $\bar \bv =\{\bar  v_n , n\in  \N\}$ defines   an infinite spine or branch. 
We denote  by $\T_1$ the subset of trees with only one infinite spine:
\begin{equation}
   \label{eq:defT1}
\T_1=\{\bt\in \T; \text{ there exists a unique } v \in (\N^* )^{\N^*} \text{ s.t. } \bar \bv 
 \subset \bt \}.
\end{equation}

  For  $h\in \N$,  the
restriction function $r_{h}$ from $\T$ to $\T^{(h)}$ is defined by:
\begin{equation}
   \label{eq:rh}
\forall \bt\in\T,\ r _{h }(\bt)=\{u\in\bt;\, H(u) \le h\}
\end{equation}
that is, $r_h(\bt)$ is the sub-tree of $\bt$ obtained by cutting the
tree at height $h$.
We endow the set $\T$ with the  distance:
\[
\delta (\bt,\bt')=2^{-\sup\{h\in\N,\ r _{h }(\bt)=r_{h }(\bt')\}}.
\]
It is  easy to  check that  this distance is in fact ultra-metric, that  is,
for all $\bt,\bt',\bt''\in\T$:
\[
\delta (\bt,\bt')\le \max\bigl(\delta(\bt,\bt''),\delta(\bt'',\bt')\bigr).
\]
Therefore all the open balls are closed. Furthermore, for $\bt\in \T_0$,
the singleton $\{\bt\}$ is also equal to the  open ball centered at
$\bt$ with radius   less than $2^{-H(\bt)}$. 
Notice also that for $\bt\in
\T$ and $h\in \N$, the set: 
\begin{equation}
   \label{eq:rhrh}
r_h^{-1}(\{r_h(\bt)\})=\{\bt'\in \T;\, \delta(\bt, \bt')\leq 2^{-h}\}
\end{equation}
is the (open and closed) ball centered at $\bt$ with radius $h$. 
 The restriction functions
are contractant with respect to the distance $\delta$ and thus
continuous.

Let $u\in  \cu$.  Recall $k_u(\bt)$ is  the number of offsprings  of the
node  $u$   in  $\bt$,  with   the  convention  that   $k_u(\bt)=-1$  if
$u\not  \in \bt$.   For  $u\in  \cu$ and  $\bt,  \bt'\in  \T$ such  that
$\delta(\bt, \bt')< 2^{-H(u)}$, we get that $k_u(\bt) =k_u(\bt')$.  This
implies that the function $\bt \mapsto k_u(\bt)$ is continuous on $\T$.


A  sequence $(\bt_n,  n\in\N)$  of trees  in $\T$  converges  to a  tree
$\bt\in \T$ with respect  to the distance $\delta$ if and  only if, for every
$h\in \N$, we have $r_{h} (\bt_n)=r_{h} (\bt)$ for $n$ large
enough, and thus  if and  only if for   all   $u\in   \cu$: 
\[
\lim_{n\rightarrow+\infty   }
k_u(\bt_n)=k_u(\bt)\in  \N\cup  \{-1   \}. 
\]

We end this section by stating that $\T$ is a Polish metric space (but
not  compact), that is a complete separable metric space. 
\begin{lmm}
   \label{lem:polish}
The metric space  $(\T ,\delta )$  is a
Polish metric space.
\end{lmm}

\begin{proof}
Notice that  $\T_0$, which is countable,  is dense in $\T $ as for all
$\bt\in \T$, the sequence  $(r_h(\bt), h\in \N)$ of elements of $\T_0$
converges to $\bt$. So the metric space $(\T ,\delta )$ is separable. 

Let $(\bt_n, n\in \N)$ be a Cauchy sequence in $\T$. Then for all $h\in
\N$, the sequence $(r_h(\bt_n), n\in \N)$ is a Cauchy sequence in
$\T^{(h)}$. Since for $\bt, \bt'\in \T^{(h)}$, $\delta(\bt, \bt')\leq
2^{-h}$ implies that $\bt=\bt'$, we deduce that the sequence
$(r_h(\bt_n), n\in \N)$ is constant for $n$ large enough equal to say
$\bt^{h}$. By continuity of the restriction functions, we deduce that
$r_h(\bt^{h'})=\bt^h$ for any $h'>h$. This implies that $\bt=\bigcup _{h\in
  \N} \bt ^h$ is a tree such that $r_h(\bt)=\bt^h$ for all $h\in\N$, and that the sequence $(\bt_n, n\in \N)$
converges to $\bt$. This gives that the metric space  $(\T ,\delta)$  is complete. 
\end{proof}

\subsection{Bienaym\'e-Galton-Watson trees}
\label{sec:GW}
\subsubsection{Definition}
\label{sec:def}
Let $p=(p(n), n\in \N)$ be a  probability distribution on the set of the
non-negative integers and $\zeta$ be a random variable with distribution
$p$.    Let  $g_p(r)=\E\left[r^\zeta\right]$,   $r\in  [0,1]$,   be  the
generating  function  of  $p$,  $\rho(p)$  its  convergence  radius  and
$m(p)=g_p'(1)=\E[\zeta] $ its mean which belongs to $[0, +\infty ]$.  We
will write $g$,  $\rho$ and $m$ for $g_p$, $\rho(p)$  and $m(p)$ when it
is clear from the context. Let $A\subset \N$ be not empty nor reduced to
$\{0\}$, and  write $\text{GCD}(A)$ for  the greatest common  divisor of
the             integers             in            $A$.          The period 
of $p$ is defined by:
\[
d=\max\{k;  \,  \supp(p)\subset  k\N\}=\text{GCD}  (\supp(p)).
\]
We say that $p$ is aperiodic if $d=1$.

\begin{dfntn}[Branching property and BGW tree]
   \label{defi:GW-T}
A $\T$-valued random variable $\tau$ is 
said to satisfy the branching property if  for $n\in  \N^*$, conditionally  on 
$\{k_\root(\tau)=n\}$,                 the                 sub-trees
$(\cs_1(\tau),\cs_2(\tau),\ldots,\cs_n(\tau))$   are   independent   and
distributed as the original tree $\tau$. 

A $\T$-valued random variable $\tau$ is a BGW tree
with   offspring distribution     $p$   if  it satisfies the branching
property and   the    distribution   of
$k_\root(\tau)$  is $p$. 
\end{dfntn}

It is easy to check that $\tau$ is a BGW tree with offspring distribution $p$
if and only if  for every  $h\in \N^*$ and  $\bt\in
\T^{(h)}$, we have:
\begin{equation}
\label{eq:loi-tau_h}
\P(r_{h }(\tau)=\bt) =\prod_{u\in \bt, \, \val{u} <h }
p\bigl(k_u(\bt)\bigr). 
\end{equation}
In particular, the restriction of the distribution of $\tau$ on the
set $\T_0$ is given by: 
\begin{equation}
\label{eq:loi-tau}
\forall \bt\in \T_0,\quad \P(\tau=\bt)=\prod_{u\in\bt}p\bigl(k_u(\bt)\bigr).
\end{equation}

It is easy to check the following lemma. Recall the definition of the
BGW process $Z=(Z_h, h\in \N)$ given in \reff{eq:def-Z}. 
\begin{lmm}[BGW process]
   \label{lem:z=Z}
Let $\tau$ be a BGW tree. The process $(z_h(\tau), h\in \N)$ is
distributed as $Z$.  
\end{lmm}

The offspring distribution $p$ and the BGW  tree are called critical (resp.  sub-critical, super-critical) if
$m(p)=1$ (resp. $m(p)<1$, $m(p)>1$). 

\subsubsection{Extinction probability}
\label{sec:ce}
Let $\tau$ be a BGW tree  with offspring distribution $p$. The extinction
event of  the BGW tree  $\tau$ is  $\ce(\tau)=\{\tau\in \T_0\}$, which we
shall denote $\ce$ when there is no possible confusion. Thanks to
Lemma \ref{lem:z=Z},  this is coherent with  Definition~\eqref{eq:def-E}.
We have the following particular cases:
\begin{itemize}
  \item If $p(0)=0$, then $\P(\ce)=0$ and a.s.\ $\tau\not\in \T_0$. 
    \item If $p(0)=1$, then a.s.\ $\tau=\{\root\}$ and $\P(\ce)=1$. 
    \item If    $p(1)=1$, then $m(p)=1$  and a.s.\  the tree
      $\tau=\bigcup  _{n\ge  0}{\{1\}^n}$,   with  the  convention  that
      $\{1\}^0= \{\root\}$, is  reduced to one infinite  spine.  In this
      case $\P(\ce)=0$.
\item If $0<p(0)<1$ and $p(0)+p(1)=1$, then  $H(\tau)+1$ is a geometric random
     variable with parameter $p(0)$ and $\tau=\bigcup
     _{0\leq n\leq H(\tau)}{\{1\}^n}$. In this case
     $\P(\ce)=1$. 
\end{itemize}
From now  on, we shall  omit the  particular critical case $p(1)=1$.
Let  $q$ be  the smallest root of the equation $g(r)=r$ in $[0, 1]$
(which exists since $g(1)=1$ and $g$ is continuous). 

\begin{figure}[H]
\begin{center}
\psfrag{G(0)}{\hspace{-.4cm} $p(0)$}
\psfrag{0}{}
\psfrag{00}{$0$}
\psfrag{1}{}
\psfrag{11}{$1$}
\psfrag{q}{$q$}
\hfill\includegraphics[height=4cm]{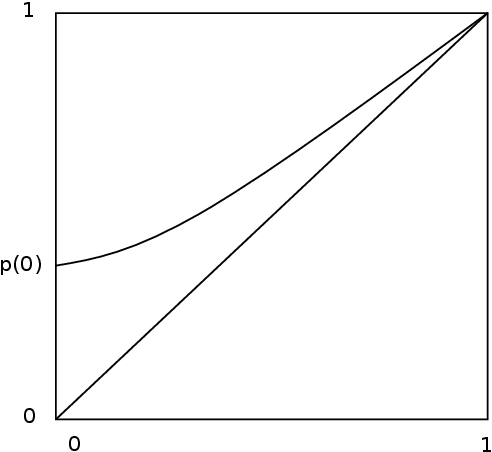}
\hfill\includegraphics[height=4cm]{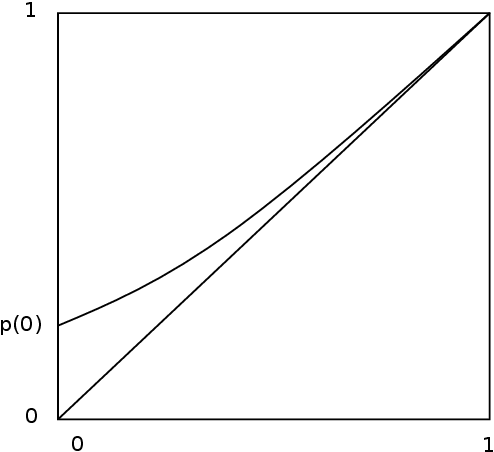}
\hfill\includegraphics[height=4cm]{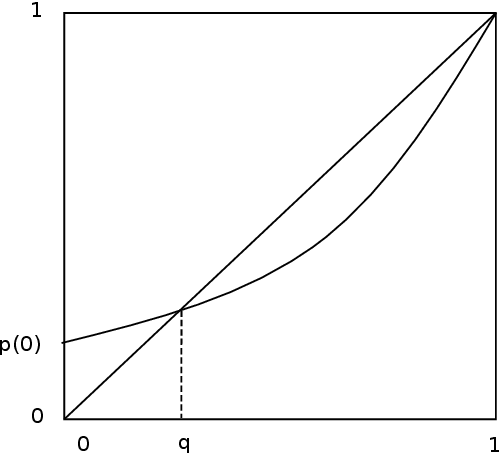}
\end{center}
\caption{Generating function in the sub-critical (left), critical
  (middle) and super-critical (right) cases.}
\label{fig:gene1}
\end{figure}

\begin{rmrk}[Roots of $g(r)=r$]
   \label{rem:prop-g}
We  clearly get that $q=0$ if and only if  $p(0)=0$. 
We now suppose that $p(0)>0$. If $p(0)+p(1)=1$, then we clearly get $q=1$. 

If $p(0)+ p(1)<1$,  we get that $g$ is strictly convex  the equation
$g(r)=r$ has at most two roots in $[0,1]$. 
 Since $g(1)=1$, we get  that  $q=1$ if $g'(1)\leq 1$ 
(in this case the equation $g(r)=r$ has only one root in $[0,1]$) and,
since furthermore $g(0)\geq 0$, we also get that 
$0\leq q<1$ if $g'(1)>1$, see Fig.~\ref{fig:gene1}.
\end{rmrk}

Using the branching property,  we get:
\begin{align*}
   \P(\ce) & =\P(\ce(\tau))\\
& =\sum_{k \in \N} \P\bigl(\ce\bigl(\cs_1(\tau)\bigr), \ldots, \ce
\bigl(\cs_k(\tau)\bigr)\bigm| 
k_\root(\tau)=k\bigr) p(k)\\
& = \sum_{k \in \N} \P(\ce)^k p(k)\\
& =g\bigl(\P(\ce)\bigr).  
\end{align*}
We deduce that $\P(\ce)$ is a root in $[0,1]$ of the equation
$g(r)=r$. The following corollary is then an immediate consequence of
Remark  \ref{rem:prop-g}.

\begin{crllr}[(Sub-)critical case]
   \label{cor:sub-critical}
   When the offspring distribution $p$ is critical with $p(1)<1$ or
   sub-critical, we have a.s.\ extinction for the associated BGW tree, that is, $\P(\ce)=1$. 
\end{crllr}

Let   $p$  be   a  super-critical   offspring  distribution   satisfying
$p(0)>0$. In this case we have $0<q<1$, and thus $\P(\ce)>0$. For
$n\in \N$, we set:
\begin{equation}
   \label{eq:tilde-p}
\tilde p(n)=q^{n-1}p(n).
\end{equation}
Since $\sum_{n\in \N} \tilde  p(n)={g(q)}/{q}=1$, we deduce that $\tilde
p=(\tilde p(n),  n\in \N)$ is  a probability distribution on  $\N$. Since
$g_{\tilde   p}   (r)=g(qr)/q$,   we    deduce   that   $g_{\tilde   p}'
(1)=g'(q)<1$. This implies that the offspring distribution $\tilde p$ is
sub-critical.   In
particular, if  $\tilde \tau$ is  a BGW tree with  offspring distribution
$\tilde p$, we have $\P(\ce(\tilde \tau))=1$.

\begin{lmm}[BGW tree conditioned on extinction]
   \label{lem:sur-critique}
   For a super-critical BGW tree $\tau$ with offspring distribution $p$,
   we  have $\P(\ce)=q$.  Furthermore, if $p(0)>0$, then conditionally  on the  extinction
   event, $\tau$ is distributed as a sub-critical BGW tree $\tilde \tau$
   with offspring distribution $\tilde p$ given by \reff{eq:tilde-p}.
\end{lmm}

\begin{proof}
  If $p(0)=0$, then we have $\P(\ce)=0$ and $q=0$.
  We now assume that $p(0)>0$. According to Corollary \ref{cor:sub-critical}, a.s.\ $\tilde \tau$
belongs to $\T_0$. For $\bt\in \T_0$, we have:
\[
   q\P(\tilde \tau=\bt)
= q\prod_{u\in \bt} p(k_u(\bt)) q^{k_u(\bt) -1} 
= q^{1+\sum_{u\in \bt}(k_u(\bt) -1)} \prod_{u\in \bt} p(k_u(\bt))
= \P( \tau=\bt),
\]
where we used \reff{eq:loi-tau} and the definition of $\tilde p$ for the
first equality and \reff{eq:sum_k} as well as \reff{eq:loi-tau} for the
last one. We deduce, by summing the previous equality over all finite trees $\bt\in\T_0$ that, for any non-negative function $\ch$ defined on $\T_0$,
\[
\E\left[\ch(\tau) \ind_{\{\tau\in \T_0\}}\right]
= q\E\left[\ch(\tilde \tau) \right],
\]
as $\tilde \tau$ is a.s. finite. Taking $\ch=1$, we deduce that
$\P(\ce(\tau))=q$. Then we get:
\[
\E\left[\ch(\tau) |\ce(\tau) \right]
= \E\left[\ch(\tilde \tau) \right].
\]
Thus, conditionally on the
extinction event,  $\tau$  is distributed as $\tilde \tau$. 
\end{proof}

We deduce the following corollary on  BGW processes. 
\begin{crllr}[BGW process conditioned on extinction]
   \label{cor:cond-Z-sur}
Let $Z$ be  a super-critical  BGW process  with offspring distribution $p$
satisfying $p(0)>0$. Conditionally on the
extinction event,  $Z$ is distributed as a sub-critical BGW process
$\tilde Z$ with offspring distribution $\tilde p$ given by~\reff{eq:tilde-p}.  
\end{crllr}

\subsubsection{Distribution of the super-critical BGW tree conditionally
   on the  non-extinction event}
\label{sec:super-cec}

Let $\tau$ be a super-critical BGW tree with offspring distribution $p$.
We  shall  present  a  decomposition  of  the  super-critical  BGW  tree
conditionally on the  non-extinction event $\ce^c=\{H(\tau)=+\infty \}$.
Notice that the  event $\ce^c$ has positive probability  $1-q$, with $q$
the smallest root of $g(r)=r$ on $[0,1]$.

We say that $v\in \bt$ is a \textit{survivor} in $\bt\in \T$ if
$\Card(S_v(\bt))=+\infty $ and becomes \textit{extinct}
otherwise. We define the survivor process $(z_h^s(\bt), h\in \N)$ by:
\[
z_h^s(\bt)=\Card \bigl(\{u\in \bt; \, H(u)=h \text{ and $u$ is a survivor}\}\bigr).
\]

Notice that the root $\root$ of $\tau$ is
a survivor
with probability $1-q$. Let $S$ and $E$ denote respectively the numbers of
children of the root which are survivors and which become extinct. We define for $r, \ell\in [0,1]$:
\[
G(r,\ell)=\E\left[r^S \ell^E|\ce^c\right].
\]

We have the following lemma. 

\begin{lmm}
   \label{lem:NE}
   Let $\tau$  be a super-critical  BGW tree with  offspring distribution
   $p$  and let $q$ be the smallest root of $g(r)=r$
   on $[0,1]$. We have for $r, \ell\in [0,1]$:
\begin{equation}
   \label{eq:def-G}
G(r,\ell)=\frac{g\bigl((1-q)r+q\ell\bigr) - g(q\ell)}{1-q}\cdot
\end{equation}
\end{lmm}
\begin{proof}
Recall $q<1$.    We have:
\begin{align*}
   \E\left[r^S \ell^E|\ce^c\right]
&=\inv{1-q} \E\left[r^S \ell^E\ind_{\{S\geq 1\}}\right]\\
&=\inv{1-q} \sum_{n\in \N^*} p(n) \sum_{k=1}^n \binom{n}{k} (1-q)^k r^k
q^{n-k} \ell^{n-k} \\ 
&=\inv{1-q} \sum_{n\in \N^*} p(n) \Bigl(\bigl((1-q)r + q\ell\bigr)^n - (q
  \ell)^n\Bigr)\\
&=\frac{g\bigl((1-q)r+q\ell\bigr) - g(q\ell)}{1-q},
\end{align*}
where we used the branching property and the fact that a BGW tree with
offspring distribution $p$ is finite with probability $q$ in the second equality.
\end{proof}

We consider the following two-type BGW tree $\mt^s$ distributed as
follows:

\begin{itemize}
   \item[-] Individuals are of type
$s$ (for \emph{survivor})  or of type $e$ (for \emph{extinct}).
   \item[-] The root of $\mt^s$ is of type $s$. 
   \item[-] An individual of type $e$ produces only individuals of type
     $e$ according to the sub-critical offspring distribution $\tilde p$ defined by
     \reff{eq:tilde-p}. 
\item[-] An individual of type $s$ produces $S\geq 1$ individuals of type $s$ and
  $E$ of type $e$, with generating function $\E\left[r^S
    \ell^E\right]=G(r,\ell)$ given  by \reff{eq:def-G}. 
Furthermore the order of the $S$ individuals of type $s$ and of
the $E$ individuals of type $e$ is uniform among the $\binom{E+S}{S}$ possible
configurations. Thus the probability for an individual $u$ of type $s$ to
have $n$ children and  whose children of type $s$ are  $\{ui, \, i\in A\}$, with
$A$ a non-empty subset $\{1, \ldots, n\}$ of cardinal $\sharp A$, is
(with the convention $0^0=1$):
\begin{equation}
   \label{eq:pnA}
p(n) (1-q)^{\sharp A-1} q^{n-\sharp A}.
\end{equation}
This indeed define a probability measure as:
\begin{align*}
\sum_{n\in \N^*} \sum_{A\subset \{1, \ldots, n\}, \, A\neq
  \emptyset}  p(n) (1-q)^{\sharp A-1} q^{n-\sharp A}
&= \sum_{n\in \N^*} \sum_{k=1}^n p(n)\binom{n}{k} (1-q)^{k-1} q ^{n-k}\\
&=  \sum_{n\in \N^*} p(n) \frac{1- q^n}{1-q}\\
&= \frac{g(1) - g(q)}{1-q}=1.
\end{align*}
\end{itemize}

Notice that an individual in $\mt^s $ is a survivor if and only if it is
of type $s$. We write $\tau^s$ for the $\T$-valued random variable
defined as  $\mt^s$ when
forgetting the types. 

Using the branching property, it is easy to deduce the following
corollary. 

\begin{crllr}[BGW tree conditioned on non extinction]
   \label{cor:TS}
  Let $\tau$  be a super-critical  BGW tree with  offspring distribution
   $p$. Conditionally on $\ce^c$, $\tau$ is distributed as
   $\tau^s$. 
\end{crllr}

\begin{proof}
If $p(0)=0$, then we have $q=0$,  $\P(\ce^c)=1$ and 
$\tau^s=\tau$.

We now assume that $q>0$.
We denote by $\cs_h$ the set of individuals  of $\tau^ s$ 
    at height $h$ whose type in $\mt^s$ is  $s$.  Because ancestors of an  individual of type $s$ are also of
  type $s$ and that every individual of type $s$ has at least a child of
  type  $s$, we  deduce that  $\mt^s$  truncated at  level $h$  is
  characterized by $r_h(\tau^s)$ and $\cs_h$.

Let $\bt\in \T_0$ such that $H(\bt)=h$ and $A\subset \{u\in \bt, \, H(u)=h\}$ with $A\neq
\emptyset$. Set $n=z_h(\bt)$. Let  $\ca=\bigcup _{u\in A} \{v\in \cu; \,
  v\prec u\}$ be the set of ancestors of elements of 
  $A$
and set $\ca^c=r_{h-1}(\bt)\backslash \ca$. 
For $u\in \ca$, we denote by $k_u^s(\bt, A)$ the number of children of
$u$ in $\bt$ that belong
to $\ca\cup  A$. 
We have:
\begin{align*}
   \P(r_h(\tau^s)  =\bt, \, \cs_h=A)
&=\prod_{u\in \ca^c} \tilde p (k_u(\bt)) \prod_{u\in \ca} p(k_u(\bt))\, 
  (1-q)^{k_u^s(\bt, A) -1} q ^{k_u(\bt) - k_u^s(\bt, A)}\\
&=\left(\prod_{u\in \ca^c\cup \ca} p (k_u(\bt)) \right)
\, \left(\frac{1-q}{q}\right)^{\sum_{u\in \ca} (k_u^s (\bt, A)-1)}
\, q^{\sum_{u\in \ca\cup \ca^c} (k_u (\bt)-1)}\\
&=\P(r_h(\tau)=\bt) \, \left(\frac{1-q}{q}\right)^{\sharp A-1} \, q^{n -1}\\
&=\inv{1-q}\, \P(r_h(\tau)=\bt) \, (1-q)^{\sharp A} q^{n-\sharp A}, 
\end{align*}
where we  used \reff{eq:pnA} for the  first equality, Definition \eqref{eq:tilde-p} of $\tilde p$ for the second one and, for  the third
equality,  Formula  \reff{eq:sum_k}  twice   as  well  as  $n=z_h(\bt)$.
Summing the previous equality over all  possible choices for $A$, we get
(recall that $A$ is non empty):
\begin{align*}
   \P(r_h(\tau^s)=\bt)
&=\sum_{A} \inv{1-q}\, \P(r_h(\tau)=\bt) \, (1-q)^{\sharp A }
                         q^{n-\sharp A}\\
&= \inv{1-q}\, \P(r_h(\tau)=\bt) \sum_{k=1}^n \binom{n}{k} \, (1-q)^{k}
  q^{n-k}\\
&= \P(r_h(\tau)=\bt)\, \frac{1- q^n}{1-q}\cdot
\end{align*}
On the other hand, we have:
\begin{multline*}
\P(r_h(\tau)=\bt|\, \ce^c)
= \frac{\P(r_h(\tau)=\bt) - \P(r_h(\tau)=\bt, \ce)}{1- \P(\ce)}\\
= \frac{\P(r_h(\tau)=\bt) - \P(r_h(\tau)=\bt)\, q^n}{1- q}
= \P(r_h(\tau)=\bt)\, \frac{1- q^n}{1-q}, 
\end{multline*}
where we used the branching property at height $h$ for $\tau$ for the
second equality. Thus we
have obtained that $\P(r_h(\tau^s)=\bt)=\P(r_h(\tau)=\bt|\, \ce^c)$ for
all $\bt\in \T_0$, which concludes the proof. 
\end{proof}

In particular, it is easy to deduce from the definition of $\tau^s$,
that the survivor  process $(z_h^s(\tau), h\in \N)$ is conditionally on
$\ce^c$ a BGW process whose offspring distribution $ \hat p$ has generating function:
\[
g_{\hat p}(r)=G(r,1)=\frac{g((1-q)r+q) - q}{1-q}\cdot
\]
The mean of $ \hat p$ is  $g_{\hat p}'(1)=g'(1)$ the mean of $p$. Notice
also  that  $\hat  p(0)=0$,  so  that  the  BGW  process  with  offspring
distribution  $\hat  p$  is  super-critical and  a.s.\  does  not  suffer
extinction.

\begin{figure}[H]
\begin{center}
\psfrag{0}{}
\psfrag{00}{$0$}
\psfrag{1}{}
\psfrag{11}{$1$}
\psfrag{q}{$q$}
\psfrag{G(0)}{$p(0)$}
\includegraphics[height=5cm]{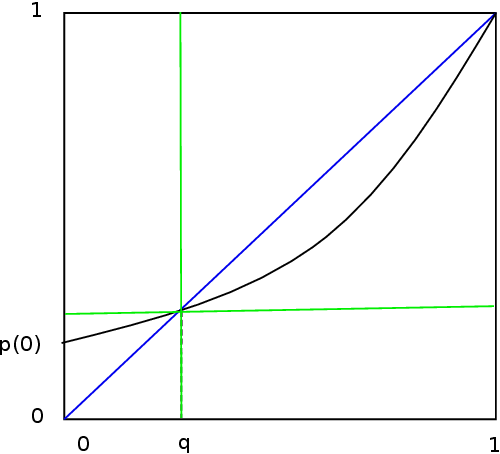}
\end{center}
\caption{In the super-critical case with $p(0)>0$ (and thus $0<q<1$): the generating functions $g_p$, $g_{\tilde
    p}$  in the lower sub-square and $g_{\hat
    p}$ in the upper sub-square (up to a scaling factor).}
\label{fig:gene2}
\end{figure}

If $p(0)>0$ (and thus $q>0)$, recall that the BGW tree $\tau$ conditionally on
the extinction event is a BGW tree with  offspring  distribution
$\tilde p$, whose generating function is:
\[
g_{\tilde p}(r)=\frac{g(qr)}{q}\cdot
\]
We observe that the generating function $g(r)$ of the super-critical
offspring distribution $p$ can be 
recovered from the extinction probability $q$, the generating functions 
$g_{\tilde p}$ 
and 
$g_{\hat   p}$ 
of the offspring distribution of
the BGW tree conditionally on the extinction event (for $r\leq q$)
 and of
the backbone process (for $r\geq q$):
\[
g(r)=qg_{\tilde p} \left(\frac{r}{q} \right)\ind_{[0,q]} (r) + 
\left(q+ (1-q) g_{\hat 
  p} \left(\frac{r-q}{1-q} \right)\right) \ind_{(q,1]} (r).
\]
We  can  therefore read  from  the  super-critical generating  functions
$g_p$,  the  sub-critical generating  function  $g_{\tilde  p}$ and  the
super-critical   generating   function    $g_{\hat   p}$,   see   Fig.~\ref{fig:gene2}.

\subsection{Kesten's tree}
\subsubsection{Definition}
\label{sec:defK}
The  Kesten's tree  is a multi-type
Galton-Watson tree, that is a random tree where all individuals
reproduce independently of the others, but the offspring distrbution
depends on the type of the individual.
For a probability distribution $p=(p(n), n\in \N)$ on $\N$ with finite
positive mean $m$, the
corresponding size-biased distribution $p^*=(p^*(n), n\in \N)$ is
defined by:
\[
p^*(n)=\frac{np(n)}{m}\cdot
\]
The two-type BGW Kesten's tree $\mt^*$  associated with the
probability distribution $p$ on $\N$ 
is  distributed as follows: 
\begin{itemize}
   \item[-] Individuals are \textit{normal} or \textit{special}.
   \item[-] The root of $\mt^*$ is   \textit{special}.
   \item[-] A \textit{normal} individual   produces only
     \textit{normal} individuals according 
     to  $p$. 
\item[-] A \textit{special}  individual  produces individuals 
  according to the  size-biased distribution $p^*$ (notice that it
   has always at least one offspring since $p^*(0)=0$). One of them,
  chosen uniformly at random, is  
  \textit{special},  the others (if any) are \textit{normal}. 
\end{itemize}
\begin{dfntn}[Kesten's tree]
   \label{defi:k-tree}
Let $p$ be an offspring distribution 
 with finite positive mean ($m\in (0, +\infty) $).
The Kesten's tree $\tau^*$ associated with $p$ is  the $\T$-valued random variable
defined as  $\mt^*$ when
forgetting the types.
\end{dfntn}

Notice $\tau^*$ belongs  a.s.\ to $\T_1$, the set of  tree with only one
infinite spine, if $p$ is sub-critical or critical. In the next lemma we
provide a link between the distribution of $\tau$ and of $\tau^*$.

\begin{lmm}[Relation between the BGW tree and the Kesten's tree]
   \label{lem:distrib-ktree}
Let  $p$ be an offspring distribution 
 with finite positive mean, $\tau$
be a
BGW tree with offspring distribution $p$ and  $\tau^*$ be a Kesten's tree
associated with $p$. For all  $n\in \N$,  $\bt\in \T_0$  and $v\in \bt$ such that
$H(\bt)=H(v)=n$, we have: 
\begin{align}
\nonumber
\P\left(r_n(\tau^*)=\bt, \text{$v$ is \textit{special}}\right)
&=\inv{m^n}\P\left(r_n(\tau)=\bt\right),\\
   \label{eq:t*t}
\P\left(r_n(\tau^*)=\bt\right)
&=\frac{z_n(\bt)}{m^n} \P\left(r_n(\tau)=\bt\right).
\end{align}
\end{lmm}

\begin{proof}
   Notice  that if $u$ is special, then the probability that it has
   $k_u$ children and  $ui$ is 
   special (with $i$ given and $1\leq i\leq k_u$) is just $p^*(k_u)k_u^{-1}=p(k_u)/m$. 
Let  $n\in \N$,  $\bt\in \T^{(n)}$  and $v\in \bt$ such that
$H(\bt)=H(v)=n$. Using \reff{eq:loi-tau_h}, we have: 
\[
\P\left(r_n(\tau^*)=\bt, \text{ $v$ is special}\right)
=\prod_{u\in \bt \setminus \rA_v, \, H(u)<n} p(k_u(\bt))\prod_{u\in  \rA_v} \frac{p(k_u(\bt))}{m}
=\inv{m^n} \P\left(r_n(\tau)=\bt\right).
\]
Since there is only one special element of $\bt$ at level $n$ among the
$z_n(\bt)$ elements of $\bt$ at level $n$, summing the previous equality
over $v$ gives  \reff{eq:t*t}. 
\end{proof}

We suppose that $m\in (0, +\infty )$. We consider the  filtration $\cf=(\cf_n, n\in \N)$  generated by $\tau$:
$\cf_n=\sigma(r_n(\tau))$    and    the     normalized    BGW    process
$(W_n, n\in \N)$ defined by:
\begin{equation}
   \label{eq:defW}
W_n=\frac{z_n(\tau)}{m^n}\cdot
\end{equation}
Notice that $W_0=1$.  If necessary, we shall write $W_n(\bt )=z_n(\bt)/m^n$ to stress
the dependence in $\bt$. 

\begin{crllr}[The martingale $(W_n)$]
   \label{cor:Wmart}
   Let $p$ be an  offspring distribution  with
   finite  positive mean.  The process  $(W_n, n\in \N)$  is  a  non-negative
   martingale adapted to the filtration~$\cf$.
\end{crllr}

\begin{proof}
   Let $\rP$ and $\rP^*$ denote respectively the distribution on $\T$ of a BGW
tree $\tau$ with offspring distribution $p$ and Kesten's tree $\tau^*$
associated with $p$. We deduce from Lemma
\ref{lem:distrib-ktree} that for all $n\geq 0$:
\begin{equation}
   \label{eq:P*P}
d\rP^*_{|\cf_n}(\bt)=W_n(\bt)\,  d\rP_{|\cf_n}(\bt).
\end{equation} 
This implies that $(W_n, n\geq 0)$ is a non-negative $\rP$-martingale adapted
to the filtration $\cf$. 
\end{proof}

\subsubsection{Asymptotics  of the BGW process}
\label{sec:asymptZ}

Let $\tau$ be  a  BGW tree with  offspring distribution $p$
 having finite positive mean $m=g'(1)$ and $p(1)<1$.  Recall
the renormalized BGW process $(W_n,  n\in \N)$ defined by \reff{eq:defW}
and $q$ the smallest root of $g(r)=r$ in $[0, 1]$.

\begin{lmm}[On $W =\lim_{n\rightarrow \infty } W_n$]
   \label{lem:a-Pw=0}
Let $p$ be an  offspring distribution 
with finite positive  mean.
The sequence $(W_n, n\in
\N)$ converges a.s.\ to a random variable $W $ such that $\E[W ]\leq 1$ and
$\P(W =0)\in \{q, 1\}$. 
\end{lmm}

\begin{proof}
According to Corollary \ref{cor:Wmart},  $(W_n, n\in  \N)$ is  a non-negative
martingale.  Thanks  to the  convergence  theorem  for  martingales,  see
Theorem 4.2.10  in \cite{d:pte}, we  get that it converges a.s.\ to  a non-negative
random variable $W $ such that $\E[W ]\leq 1$.

By decomposing $\tau$ with respect to the children of the root, we get:
\[
W_n(\tau) =\inv{m} \sum_{i=1}^{k_\root(\tau)} W_{n-1} (\cs_i(\tau)).
\]
The    branching     property    implies    that     conditionally    on
$k_\root(\tau)$,  the  random   trees  $\cs_i(\tau)$,  $1\leq  i\leq
k_\root(\tau)$  are  independent  and  distributed  as  $\tau$.   In
particular, $(W_n(\cs_i(\tau)), n\in \N)$ converges a.s.\ to a limit, say
$W^i $,  where $(W^i ,  i\in  \N^*)$ are  independent non-negative  random
variables distributed  as $W $ and independent  of $k_\root(\tau)$. By
taking the limit as $n$ goes to infinity, we deduce that a.s.:
\[
W =\inv{m} \sum_{i=1}^{k_\root(\tau)} W ^i. 
\]
This implies that:
\[
\P(W =0)=\sum_{n\in \N} p(n) \P(W ^1=0, \ldots, W ^n=0)=g(\P(W =0)).
\]
This implies that $\P(W =0)$ is a non-negative solution of $g(r)=r$ and
so belongs to $\{q, 1\}$. 
\end{proof}

\begin{rmrk}[The (sub-)critical case]
   When the offspring distribution  $p$ is critical with $p(1)<1$ or  sub-critical, we have $q=0$, 
 a.s.\   $W_n=0$ for $n$ large enough, and thus a.s.\ $W=0$. 
\end{rmrk}

We aim now to compute $\P(W =0)$ in the super-critical case. The following result goes back to Kesten
and Stigum \cite{ks:ltmgwp} and we present the proof by Lyons, Pemantle and
Peres \cite{lpp:cp}. 
Recall  that $\zeta$ is  a random variable with distribution $p$. We use the
notation  $\log^+(r)=\max(\log(r), 0)$.

\begin{thrm}[The $L \log L$ condition]
   \label{theo:ks}
Let $p$ be a super-critical offspring distribution 
 with finite mean ($m\in (1,+\infty) $).
Then we have:
\[
  \P(W =0)=q
  \quad\text{if and only if}\quad
  \E\left[\zeta \log^+(\zeta)\right] <+\infty.
\]
\end{thrm}
In particular, we also have that $\P(W =0)=1$ if and only if
$\E\left[\zeta \log^+(\zeta)\right]$ is infinite.

\begin{rmrk}[Exponential growth of the population in the super-critical case]
   \label{rem:W=0}
   Assume  that  $p$   is  super-critical  and  satisfy the $L \log L$
   condition.  Since
   $\ce\subset \{W =0\}$ and $\P(\ce)=q$, we deduce that on the survival
   event  $\ce^c$  a.s.\  $\lim_{n\rightarrow+\infty  }  z_n(\tau)/m^n=W
   >0$. (On the extinction event $\ce$, we have that a.s.\ $z_n(\tau)=0$
   for $n$ large.)  So, a.s.\ on the survival event, the population size
   at level  $n$ behaves  like a positive  finite random  constant times
   $m^n$.
\end{rmrk}

\begin{proof}[Proof of Theorem~\ref{theo:ks}]
We use notations from the proof of Corollary \ref{cor:Wmart}:  $\rP$ and $\rP^*$ denote respectively the distribution of a BGW
tree $\tau$ with offspring distribution $p$ and of a Kesten's tree $\tau^*$
associated with $p$. According to \reff{eq:P*P}, for all $n\geq 0$:
\[ 
d\rP^*_{|\cf_n}=W_n\,  d\rP_{|\cf_n}.
\]
This implies that $(W_n, n\geq 0)$ converges $\rP$-a.s.\ (this is already
in  Lemma \ref{lem:a-Pw=0})  and $\rP^*$-a.s.\  to $W$  taking values  in
$[0,+\infty ]$. According  to Theorem 4.3.3 in \cite{d:pte}, we get that
for any  measurable subset $B$
of $\T$:
\[
\rP^*(B)= \rE[W\ind_B]+ \rP^*(B, W=+\infty ).
\]
Taking $B=\Omega$ in the previous equality gives:
\begin{equation}
   \label{eq:PPA}
\begin{cases} \rE[W]=1 \Leftrightarrow \rP^*(W=+\infty )=0 \quad\text{and}\\
\rP(W=0)=1 \Leftrightarrow \rP^*(W=+\infty )=1.
\end{cases}
\end{equation}
So we shall study the behavior of $W$ under $\rP^*$, which turns out to
be (almost) elementary. We first use a similar description as
\reff{eq:def-Z} to describe $(z_n(\tau^*), n\in \N)$. 

Recall that $\zeta$  is a random variable with  distribution $p$. Notice
that $p^*(0)=0$,  and let $Y$ be  a random variable such  that $Y+1$ has
distribution  $p^*$.  Under $\P$, let  $(\zeta_{i,  n};   i\in  \N,  n\in  \N)$  be
independent  random variables  distributed as  $\zeta$, and  let  $(Y_n, n\in
\N^*)$  be   independent  random   variables  distributed  as   $Y$  and
independent of $(\zeta_{i, n}; i\in \N,  n\in \N)$. We set $Z^*_0=0$ and
for $n\in \N^*$:
\[
Z^*_{n}=Y_n+\sum_{i=1}^{Z^*_{n-1}} \zeta_{i, n},
\]
with the convention that $\sum_\emptyset=0$. In particular $(Z_n^*, n\in
\N)$ is under $\P$ a BGW process with immigration 
with offspring distribution $p$ and immigration
distributed as $Y$. By construction $(Z_n^*+1, n\in \N)$ under $\P$ is distributed
as $(z_n(\tau^*), n\in \N)$ under $\rP^*$. We deduce that $(W_n, n\in \N)$ is under $\rP^*$
distributed as $(W^*_n+ m^{-n}, n\in \N)$ under $\P$, with $W^*_n=Z^*_n/m^n$. 

\medskip

Let  $(X_n,   n\in  \N)$  be  random  variables
distributed as a non-negative random variable $X$.
We recall  the following  result, which  can be  deduced from  the
Borel-Cantelli  lemma using that $\sum_{n\in \N^ *} \P(X_n/n\geq
\varepsilon)= \E[\lfloor X/\varepsilon \rfloor]$ for $\varepsilon>0$.   We have:
\begin{equation}
   \label{eq:Xn/n0}
\E[X]<+\infty \Rightarrow  \text{ a.s. } \lim_{n\rightarrow+\infty }
\frac{X_n}{n}=0. 
\end{equation}
Furthermore, if the random variables $(X_n, n\in \N)$ are
independent, then:
\begin{equation}
   \label{eq:Xn/n+}
\E[X]=+\infty \Rightarrow \text{ a.s. } \limsup_{n\rightarrow+\infty }
\frac{X_n}{n}=+\infty .
\end{equation}

\medskip
We consider the case:  $\E\left[\zeta \log^+(\zeta)\right] <+\infty
$. This implies that $\E[\log^+(Y)]<+\infty $. And according to
\reff{eq:Xn/n0}, we deduce that for $\varepsilon>0$, $\rP^*$-a.s.\  $Y_n\leq
\expp{ n\varepsilon}$ for $n$ large enough. Denote by $\cy$ the
$\sigma$-field generated by $(Y_n, n\in \N^*)$ and by $(\cf^*_n, n\in
\N)$ the filtration generated by $(W^*_n, n\in \N)$. 
Using the branching property, it is easy to
get:
\[
\E\left[W^*_n|\cy, \cf^*_{n-1}\right]
= \inv{m^n}\E\left[Z^*_n|\cy, \cf^*_{n-1}\right]
= \frac{Z^*_{n-1}}{m^{n-1}} + \frac{Y_n}{m^n}\geq W^*_{n-1}.
\]
We deduce that $(W^*_n, n\in \N)$ is a non-negative sub-martingale with
respect to $\P(\cdot|\cy)$, that is, with respect to the filtration
$(\cy \vee \cf^*_n, n\in \N)$. We
also obtain:
\[
\E\left[W^*_n|\cy\right]
=  \sum_{k=1}^n\frac{Y_k}{m^k}\leq  \sum_{k=1}^{+\infty }\frac{Y_k}{m^k}\cdot
\]
For       $\varepsilon>0$,      we       have      that       $\P$-a.s.\
$Y_n \leq  \expp{ n\varepsilon}$  for $n$ large  enough. We  deduce that
$\P$-a.s.\ we have $\sup_{n\in  \N} \E\left[W^*_n|\cy\right] <+\infty $,
that  is,   the  sub-martingale  $(W^*_n,   n\in  \N)$  is   bounded  in
$L^1(\P(\cdot|\cy))$.   Adapting    the   proof   of    Theorem   4.2.10
in~\cite{d:pte} with respect to  the probability $\P(\cdot|\cy)$, we get
that  the  non-negative  sub-martingale  $(W^*_n,  n\in  \N)$  converges
$\P(\cdot|\cy)$-a.s.\      to      a       finite      limit.      Since
$(W^*_n+  m^{-n}, n\in  \N)$ is  distributed as  $(W_n, n\in  \N)$ under
$\rP ^*$, we get that $\rP^*$-a.s.\ $W$ is finite. Use the first part of
\reff{eq:PPA} to deduce that  $\rE[W]=1$.  Since $\rP(W=0)\in \{q, 1\}$,
see Lemma \ref{lem:a-Pw=0}, we get that $\rP(W=0)=q$.

\medskip

We consider the  case: $\E\left[\zeta \log^+(\zeta)\right] =+\infty
$.  According  to   \reff{eq:Xn/n+},  we   deduce  that   for  any
$\varepsilon>0$,  a.s.\  $Y_n\geq  \expp{ n/\varepsilon}$  for infinitely
many $n$. Since $Z_n^*\geq Y_n$ and $(W^*_n+ m^{-n}, n\in \N)$ is  distributed as
$(W_n,  n\in  \N)$  under  $\rP   ^*$,  we  deduce  that  $\rP ^*$-a.s.\  $W_n\geq
\expp{n(-\log(m) +1/\varepsilon)}$  for infinitely  many $n$.  Since the
sequence $(W_n, n\in  \N)$ converges $\rP^*$-a.s.\ to $W$ taking values in
$[0,+\infty ]$,  we deduce, by
taking $\varepsilon>0$ small enough that
$\rP^*$-a.s.\  $W=+\infty $.  Use  the second  part  of \reff{eq:PPA}  to
deduce that $\rP$-a.s.\ $W=0$. 
\end{proof}

\section{Local limits of Bienaym\'e-Galton-Watson  trees}
\label{sec:loc}

There are many kinds of limits that  can be considered in order to study
large trees, among them are the local limits and the scaling limits.
 The local limits look
at the  trees up to  an arbitrary fixed  height and therefore  only sees
what happen at a finite distance from the root. One can also look at the
local  limit  near  a  node  chosen   at  random  instead  of  near  the
root. Scaling limits consider sequences  of trees where the branches are
scaled by  some factor so that  all the nodes remain  at finite distance
from the  root. These scaling limits,  which lead to the  so-called {\sl
  continuum random trees} where  the branches have infinitesimal length,
have    been    intensively    studied     in    recent    years,    see
\cite{a:crt3,d:ltcpcgwt,dlg:rt}, see also~\cite{haas2012} for scaling limits of other random trees and~\cite{kortchemski2015condensation} for condensation phenomenon. 
Other limits have also been considered,
see for example~\cite{et:cv-tree} for a convergence to dendrons with respect to
sampling procedures.

We  will focus  in this  lecture  only on  local limits  of critical  or
sub-critical BGW trees conditioned on  being large.  The most famous type
of such a conditioning is Kesten's theorem which states that critical or
sub-critical BGW trees conditioned on reaching large heights converge
to the 
Kesten's tree  which is a  (two-type) BGW  tree with a  unique infinite
spine. This result  is recalled in Theorem \ref{thm:cv-k1}.  In order to
consider other conditionings, we  shall give in Section \ref{sec:topo},
see Proposition \ref{prop:conv_dist},  an elementary characterization of
the  local  convergence  which  is  the key  ingredient  of  the  method
presented here when the limit is the Kesten's tree.

All the  conditionings we  shall consider  can be stated  in terms  of a
functional $A(\bt)$ of the tree $\bt$ and the events we condition on are
either  of the  form  $\{A(\tau)\ge n\}$  or  $\{A(\tau)=n\}$, with  $n$
large.  Usually  handling the latter  conditioning is more  delicate. In
Section \ref{sec:criteria}, we give general assumptions on $A$ so that a
{\bf critical}  BGW tree conditioned on  such an event converges  as $n$
goes to infinity, in distribution to the Kesten's tree, see our main result,
Theorem~\ref{theo:main1}.
We then apply this result  in Section \ref{sec:appli} by considering, in
the  critical case,  the following  functional:
\begin{enumerate}
   \item[\ref{sec:Kesten}]: the  height of  the tree
     (Kesten's theorem);
   \item[\ref{sec:out}]: the maximal out-degree;
   \item[\ref{sec:horton}]: the Horton-Strahler number;
   \item[\ref{sec:large-gen}]:  the largest generation;
   \item[\ref{sec:total}]:  the total progeny;
   \item[\ref{sec:L0}]:  the number of leaves;
   \item[\ref{sec:out-d}]: and  more generally  the number of  nodes with
     given out-degree.
   \end{enumerate}
   Let us stress that in the critical case with finite variance, conditioning on the total population size (for simplicity), the number of leaves is asymptotically equal to the population size times 
the probability to have zero children, and furthermore the fluctuations are Gaussian, see 
\cite{janson:fringe,th:large-GW}. It is thus not a surprise that conditioning by the total size or the number of leaves give the same local limit. 
    The conditioning on  the total progeny or the number  of leaves where
   already known, but usually under stronger hypothesis on the offspring
   distribution (higher  moments or tail conditions). 
   We stress
   out that in Theorem~\ref{theo:main1} no  further assumptions are needed in the  critical case than
   the  non-degeneracy  condition $p(1)<1$. 
   
  The  material of this section is mainly  extracted  from~\cite{ad:llcgwtisc}.

 \medskip
 
   The  {\bf sub-critical}  case is  similar with  the functional  $A=H$
   given by  the height, but it  is much more involved  otherwise and we
   only  present here  some  results  in Section  \ref{sec:sub-critical}
   without   any   proofs.    All   the   proofs   can   be   found   in
   \cite{ad:llcgwtcc}. In  the sub-critical  case, when  conditioning on
   the number of  nodes with given out-degree, two cases  may appear. In
   the so-called  generic case presented in  Section \ref{sec:gene-sub},
   the  limiting tree  is  still a  Kesten's tree  but  with a  modified
   offspring distribution.  In the
   non-generic   case,   Section  \ref{sec:non-gene},   a   condensation
   phenomenon  occurs:  intuitively  a  node that  stays  at  a  bounded
   distance from  the root has more  and more offsprings as  $n$ goes to
   infinity;  and in  the  limit,  the tree  has  a  (unique) node  with
   infinitely many offsprings.
   This phenomenon has first been  pointed out in \cite{js:cnt} and then
   in \cite{j:sgtcgwrac} when conditioning on the total progeny.  We end
   this  chapter   by  giving   in  Section   \ref{sec:gene-non-gene}  a
   characterization of generic  and non-generic offspring distributions,
   which provides  non intuitive  behavior.

\medskip
We eventually give  further extensions in  Section
   \ref{sec:extension}.

      \medskip
   
\textbf{Take-out message}: 
Most (but not all) of the usual conditionings for a critical BGW tree to be large yield the same local limit  given by the Kesten's tree. 
The picture in the sub-critical case is more complex, with possibly different local limits:  Kesten's trees with modified offspring distribution (depending on the conditioning) or random trees with condensation. This latter case might appear  in particular  when the offspring distribution has an heavy tail. 

\goodbreak

\subsection{The topology of local convergence}
\label{sec:topo}
\subsubsection{Kesten's theorem}

We  work on  the set  $\T$  of discrete  trees with  no infinite  node,
introduced   in   Section   \ref{sec:discrete}.  Recall   that   $\T_0$,
resp. $\T^{(h)}$, denotes  the subset of $\T$ of finite  trees, resp. of
trees with height less than $h$, see \reff{eq:t0-th}. 
Recall that the
restriction function $r_{h}$ from $\T$ to $\T^{(h)}$  is defined in
\reff{eq:rh}. 
When a  sequence  of random  trees $(T_n, {n\in\N})$  converges in distribution
with respect to the distance
$\delta$ (also  called the local  topology)
toward a  random tree $T$, we shall write:
\begin{equation}
   \label{eq:def-cv-d}
T_n\overset{(d)}{\longrightarrow} T.
\end{equation}
According to \reff{eq:rhrh}, the fact that all the open balls are closed,
and  the Portmanteau  theorem (see \cite{b:cpm} Theorem 2.1), we  deduce
that if \reff{eq:def-cv-d}
then:
\begin{equation}\label{eq:cv_law}
\forall h\in\N,\ \forall
\bt\in\T^{(h)},\quad  \lim_{n\rightarrow+\infty }
\P\bigl(r_h(T_n)=\bt\bigr)
= \P\bigl(r_h(T)=\bt\bigr). 
\end{equation}
Conversely,  since  $\delta$  is  an   ultra-metric,  we  get  that  the
intersection of  two balls is a  ball (possibly empty). Thus  the set of
balls  and  the empty  set  is  a  $\pi$-system.   We deduce from  Theorem  2.3  in
\cite{b:cpm}    that   \reff{eq:cv_law}  implies
\reff{eq:def-cv-d}. Thus \reff{eq:cv_law} and
\reff{eq:def-cv-d} are equivalent. 
\medskip

Convergence  in  distribution for  the  local  topology appears  in  the
following  Kesten's theorem.  Recall that  the distribution  of the Kesten's
tree is given in Definition \ref{defi:k-tree}.

\begin{thrm}[Kesten~\cite{k:sbrwrc}]
\label{thm:cv-k1}
Let $p$ be a critical   or sub-critical offspring distribution. Let $\tau$ be a BGW tree with offspring distribution $p$
and $\tau^*$ be a Kesten's tree  associated with $p$.  For every non-negative
integer  $n$,  let $\tau_n$  be  a  random  tree distributed  as  $\tau$
conditionally on $\{H(\tau)\ge n\}$. Then we have:
\[
\tau_n\overset{(d)}{\longrightarrow} \tau^*.
\]
\end{thrm}

This theorem is stated in  \cite{k:sbrwrc} with an additional second moment
condition. The proof of the theorem stated as above is due to Janson
\cite{j:sgtcgwrac}. We will give a proof of that theorem in Section
\ref{sec:Kesten} as an application of a more general result, see Theorem
\ref{theo:main1} in the critical case and Theorem
\ref{theo:main2} in the sub-critical case.

\subsubsection{A characterization of the convergence in distribution}

Recall that for a tree $\bt\in\T$, we denote by $\cl_0(\bt)$ the set
of its leaves.
Let  $\bt\in\T$ be a tree and $x\in\cl_0(\bt)$ be a leaf of the tree
$\bt$. For $\bt'\in\T$ another tree, we denote by $\bt\circledast_x\bt'$ the tree
obtained by grafting the tree $\bt'$ on the leaf $x$ of the tree $\bt$,
that is:
\[
\bt\circledast_x\bt'=\bt\cup\{xu,u\in\bt'\}.
\]
For  $\tau$ a BGW tree with offspring  distribution $p$, we deduce from~\eqref{eq:loi-tau} 
the following useful formula:
\begin{equation}
   \label{eq:Ptt'=PtPt'-v0}
\P(\tau=\bt\circledast_x\bt')=\frac{1}{p(0)}\P(\tau=\bt)\P(\tau=\bt'). 
\end{equation}
We denote by $\T(\bt,x)$ the set of all trees obtained by grafting some tree on the leaf $x\in\cl_0(\bt)$ of
$\bt\in\T$:
\[
\T(\bt,x)=\{\bt\circledast _x\bt',\ \bt'\in\T\}.
\]
Recall that the maps $\bt
\mapsto k_u(\bt)$ are continuous for all $u\in \cu$.
We deduce that  the set $\T(\bt,x)$ is closed in $\T$ as  $\bs\in \T(\bt,x)$ if and only if
$k_u(\bs)=k_u(\bt)$ for all $u\in \bt\setminus \{x\}$.  We shall see in
Lemma \ref{lem:voisinage} below 
that  it is also open.
\medskip 

Computations of the probability of BGW trees
(or Kesten's tree) to belong to such sets are very easy and lead to
simple formulas. For example, we have for $\tau$ a  BGW tree
with offspring  distribution $p$, and  all
  finite tree $\bt\in \T_0$ and leaf $x\in \cl_0(\bt)$:
\begin{equation}
\label{eq:tau-in-T}
\P(\tau=\bt)=\P\bigl(\tau\in\T(\bt,x),\
k_x(\tau)=0\bigr)
=p(0)\, \P\bigl(\tau\in\T(\bt,x)\bigr). 
\end{equation}

The next lemma is another example of the simplicity of the
formulas.  
\begin{lmm}[Another relation  between the  BGW tree  and the  Kesten's
  tree]
   \label{lem:simple-formule} 
 Let $p$ be an offspring distribution 
  with finite positive mean.
Let 
  $\tau$ be a  BGW tree with offspring  distribution $p$ and let
  $\tau^*$ be a Kesten's tree associated with $p$. Then we have, for all
  finite tree $\bt\in \T_0$ and leaf $x\in \cl_0(\bt)$:
\begin{equation}
   \label{eq:tau*-in-T}
\P\bigl(\tau^*\in \T(\bt, x)\bigr)= \inv{m^{H(x)}} \P\bigl(\tau\in \T(\bt, x)\bigr). 
\end{equation}
\end{lmm}
In the particular case of a critical offspring distribution ($m=1$), we
get for all $\bt\in \T_0$ and  $x\in \cl_0(\bt)$:
\[
\P(\tau^*\in \T(\bt, x))=  \P(\tau\in \T(\bt, x)). 
\]
However,    we   have    $\P(\tau\in   \T_0)=1$    if   $p(1)<1$,    and
$\P(\tau^*\in \T_1)=1$, with  $\T_1$ the set of trees that  have one and
only  one  infinite  spine,  see  \reff{eq:defT1}.  (Notice  that,  when
$p(1)=1$, the tree $\tau=\tau^*$ is reduced to one spine.)

\begin{proof}
  Let $\bt\in \T_0$ and $x\in \cl_0(\bt)$. If $\tau^*\in\T(\bt,x)$, then
  the node $x$ must  be a special node in $\tau^*$ as  the tree $\bt$ is
  finite whereas  the tree $\tau^*$  is a.s.\ infinite.  Therefore, using
  arguments   similar   to   those   used  in   the   proof   of   Lemma
  \ref{lem:distrib-ktree}, we have:
\[
\P\bigl(\tau^*\in\T(\bt,x)\bigr)
=\prod_{u\in \bt \setminus (\rA_x\cup \{x\})} p(k_u(\bt))\prod_{u\in  \rA_x}
\frac{p(k_u(\bt))}{m}
=\inv{m^{H(x)}} \P\left(\tau\in\T(\bt,x)\right). 
\]
\end{proof}

The following key characterization of the local convergence is proved in
Section \ref{sec:proof}.

\begin{prpstn}
  [Characterization of the local convergence in $\T_0\cup \T_1$]
  \label{prop:conv_dist}
Let $(T_n,{n\in\N})$ and $T$ be random trees taking values in the set
$\T_0\cup \T_1$. Then the sequence $(T_n,{n\ge 0})$ converges in
distribution (for the local topology) to $T$ if and only if the two
following conditions hold:
\begin{enumerate}[(i)]
\item\label{it:cv-T0}
  For every finite tree $\bt\in\T_0$, we have $\lim_{n\rightarrow
    +\infty } \P(T_n=\bt)=
  \P(T=\bt)$.
\item\label{it:cv-Ttx}
  For   every   $\bt\in\T_0$  and    leaf
  $x\in\cl_0(\bt)$, we have $\liminf_{n\rightarrow+\infty } \P\bigl(T_n\in\T(\bt,x)
  \bigr)\geq  \P\bigl(T\in\T(\bt,x)\bigr)$.

\end{enumerate}
\end{prpstn}

\subsubsection{Proof of Proposition \ref{prop:conv_dist}}
\label{sec:proof}

We denote by $\cf$ the subclass of  Borel   sets  of   $\T$:
\[
\cf=\bigl\{\{\bt\},\ \bt\in\T_0\bigr\}\cup\bigl\{\T(\bt,x),\
\bt\in\T_0,\ x\in\cl_0(\bt)\bigr\}\cup \{\emptyset\}  .
\]

\begin{lmm}
   \label{lem:pi-system}
The family $\cf$ is a $\pi$-system. 
\end{lmm}
\begin{proof}
  Recall that  a non-empty  family of  sets is a  $\pi$-system if  it is
  stable under finite intersection.  For every $\bt_1,\bt_2\in \T_0$ and
  every $x_1\in\cl_0(\bt_1)$, $x_2\in\cl_0(\bt_2)$, we have:
\begin{equation}
   \label{eq:inter-cf}
\T(\bt_1,x_1)\cap \T(\bt_2,x_2)=
\begin{cases}
  \T(\bt_1,x_1) &
  \text{if }
 \bt_1=\bt_2\circledast
_{x_2}{\bt'_1}\quad\text{and}\quad
x_1\wedge x_2=  x_2,\\
  \T(\bt_2,x_2) &
  \text{if }
 \bt_2=\bt_1\circledast
_{x_1}{\bt'_2}\quad\text{and}\quad
x_1\wedge x_2=  x_1,\\
  \{\bt_1\cup \bt_2\} &
  \text{if }
  \bt_1=\bt\circledast
_{x_2}{\bt'_1},\ \bt_2=\bt\circledast _{x_1}{\bt'_2}
\quad\text{and}\quad x_1\neq  x_2 ,\\ 
\emptyset &
\text{in the other cases.}
\end{cases}
\end{equation}
Notice that the case $(\bt_1,x_1)=(\bt_2, x_2)$ belongs to the first two
cases. An instance of the third case is represented in
Fig.~\ref{fig:mickey}.


\begin{center}
\begin{figure}[H]

\begin{tikzpicture}

\draw[thick] (0,0) .. controls (-2,1) and (-2,2) .. (-1,3) .. controls (0,3.5) .. (1,3) .. controls (2,2) and(2,1) .. cycle;
\draw (-1,3) node{$\bullet$};
\draw (1,3) node{$\bullet$};
\draw(0,1.5) node{$\bt$};
\draw(-1,2.5) node{$x_2$};
\draw(1,2.5) node{$x_1$};
\draw[thick] (-1,3) .. controls (-2,4) and (-2,5) .. (-1,6).. controls (-0.5,5) and (-0.5,4) .. cycle;
\draw(-1.2,4.5) node{$\bt'_1$};
\draw(0,-0.3) node{$\root$};
\end{tikzpicture}
\begin{tikzpicture}

\draw[thick] (0,0) .. controls (-2,1) and (-2,2) .. (-1,3) .. controls (0,3.5) .. (1,3) .. controls (2,2) and(2,1) .. cycle;
\draw (-1,3) node{$\bullet$};
\draw (1,3) node{$\bullet$};
\draw(0,1.5) node{$\bt$};
\draw(-1,2.5) node{$x_2$};
\draw(1,2.5) node{$x_1$};
\draw[thick] (1,3) .. controls (0.5,4) and (0.5,5) .. (1,6) .. controls (2,5) and (2,4) .. cycle;
\draw(1.2,4.5) node{$\bt'_2$};
\draw(0,-0.3) node{$\root$};
\end{tikzpicture}
\begin{tikzpicture}

\draw[thick] (0,0) .. controls (-2,1) and (-2,2) .. (-1,3) .. controls (0,3.5) .. (1,3) .. controls (2,2) and(2,1) .. cycle;
\draw (-1,3) node{$\bullet$};
\draw (1,3) node{$\bullet$};
\draw(0,1.5) node{$\bt$};
\draw(-1,2.5) node{$x_2$};
\draw(1,2.5) node{$x_1$};
\draw[thick] (-1,3) .. controls (-2,4) and (-2,5) .. (-1,6).. controls (-0.5,5) and (-0.5,4) .. cycle;
\draw(-1.2,4.5) node{$\bt'_1$};
\draw[thick] (1,3) .. controls (0.5,4) and (0.5,5) .. (1,6) .. controls (2,5) and (2,4) .. cycle;
\draw(1.2,4.5) node{$\bt'_2$};
\draw(0,-0.3) node{$\root$};
\end{tikzpicture}
\caption{\label{fig:mickey}Exemple of the third case in
Equation \eqref{eq:inter-cf} where $x_1\ne x_2$,
  $\bt_1=\bt\circledast 
_{x_2}{\bt'_1}$ (on the left), $\bt_2=\bt\circledast _{x_1}{\bt'_2}$ (in
the middle) and with   $\bt_1\cup\bt_2$ (on the right).}
\end{figure}
\end{center}


Thus  $\cf$  is  stable  under   finite  intersection,  and  is  thus  a
$\pi$-system.
\end{proof}

\begin{rmrk}
The third case in Equation \reff{eq:inter-cf} was ommited in the
original paper \cite{ad:llcgwtisc} where only a special case was
considered.
\end{rmrk}

We denote by $\cf'=\{A \cap (\T_0\cup\T_1)\, \colon\, A\in
   \cf\}$ the trace of $\cf$ on
   $\T_0\cup\T_1$. 

\begin{lmm}
   \label{lem:voisinage}
   All the elements of $\cf$ are open and $\cf'$ is an open neighborhood
   system in $\T_0\cup\T_1$.
\end{lmm}

\begin{proof}
  We  first  check  that  all  the  elements  of  $\cf$  are  open.  For
  $\bt\in \T$ and $\varepsilon>0$, let $B(\bt, \varepsilon)$ be the ball
  (which is open and closed) centered at $\bt$ with radius
  $\varepsilon$. If
  $\bt\in\T_0$,  we have  $\{\bt\}=B(\bt,2^{-h})$ for  every $h>H(\bt)$,
  thus  $\{\bt\}$ is  open.  Moreover, for
  some fixed $x\in\cl_0(\bt)$, for every  $\bs\in\T(\bt,x)$ , we have:
\[
B\left(\bs,2^{-H(\bt)-1}\right)\subset \T(\bt,x),
\]
which proves that $\T(\bt,x)$ is also open.
\medskip

We check that $\cf'$ is a neighborhood system: that is, since all the
elements of $\cf$ are open,  for
all $\bt\in \T_0\cup\T_1$ and 
$\varepsilon>0$,  there exists an element of $\cf$, say $A'$, which is a
subset of $B(\bt,
\varepsilon)$ and which contains $\bt$. 

If $\bt \in \T_0$, it is enough to consider $A'=\{\bt\}$. 

Let us suppose that $\bt\in\T_1$. Let $(u_n, n\in \N^*)$ be the infinite
spine  of $\bt$  so  that $\bar  u_n=u_1  \ldots u_n  \in  \bt$ for  all
$n\in \N^*$.   Let $n\in  \N^*$ such  that $2^{-n}<\varepsilon$  and set $\bt'$ the tree  obtained by cutting the spine of $\bt$ at height $n$:
$\bt'=\{v\in  \bt;  \,   \bar  u_n  \not  \in   \rA_v\}$.  
 Notice  that
$\bar  u_n\in \bt'$,  and  set $A'=\T(\bt',  \bar u_n)  $  so that  $A'$
belongs      to      the      $\pi$-system      $\cf$.       We      get
$ \bt\in  A'\subset B(\bt, \varepsilon)$.   Thus, the trace of  $\cf$ on
$\T_0 \cup \T_1$ is a neighborhood system.
\end{proof}


We are  now ready  to prove Proposition  \ref{prop:conv_dist}, following
ideas of the proof of Theorems 2.2 and 2.3 of \cite{b:cpm}.  As $\cf$ is
countable, so  is its trace  $\cf'$ on  $\T_0\cup \T_1$. We  deduce from
Lemma~\ref{lem:voisinage} that any  open set $G$ of  $\T_0\cup \T_1$ can
be written as a countable union of some elements of $\cf'$, say $(A_i,
i\in \N)$.    For any  $\varepsilon>0$, there  exists
$n_0$                             such                             that
$\P(T\in G) \leq \varepsilon+ \P(T\in \bigcup _{i\leq n_0 } A_i)$.
Without loss of generality,  we can assume that no $A_i$  is a subset of
$A_j$   for  $1\leq   i,j\leq   n_0$  and   $i\neq   j$.  According   to
\reff{eq:inter-cf}, we get that $A_i\cap A_j$ is either empty or reduced
to a  singleton.  We  then deduce  from the  inclusion-exclusion formula
that there exists $n_1\leq  n_0$, $\bt_j\in \T_0$, $x_j\in \cl_0(\bt_j)$
for    $j\leq    n_1$,    and    $n_2<\infty$,    $\bt_\ell\in    \T_0$,
$\alpha_\ell\in  \Z$  for  $\ell\leq  n_2$ such  that,  for  any  random
variable $T'$ taking values in $\T_0\bigcup \T_1$:
\[
\P\Big(T'\in \bigcup _{i\leq n_0 } A_i\Big)
=\sum_{j\leq n_1} \P(T' \in \T(\bt_j, x_j)) + \sum_{\ell\leq n_2}
\alpha_\ell \P(T'=\bt_\ell) .
\]
We deduce, assuming that  (i) and (ii) of Proposition
\ref{prop:conv_dist} hold, that:
\[
\liminf_{n\rightarrow +\infty } 
\P(T_n\in G)\geq  \liminf_{n\rightarrow +\infty } 
\P\Big(T_n\in \bigcup _{i\leq n_0 } A_i\Big)
\geq \P\Big(T\in \bigcup _{i\leq n_0 }
A_i\Big)
\geq  \P(T\in G) -\varepsilon.
\]
Since     $\varepsilon>0$    is     arbitrary,     we    deduce     that
$\liminf_{n\rightarrow +\infty } \P(T_n\in G)\geq \P(T\in G)$ for every open set $G$ of $\T$. Thanks to
the Portmanteau theorem,  see (iv) of Theorem 2.1  in \cite{b:cpm}, we
deduce that $(T_n, n\in \N)$ converges in distribution to $T$.

\subsection{A criteria for convergence toward Kesten's tree}
\label{sec:criteria}

Using the previous lemma, we can now state a general result for
convergence of conditioned BGW trees toward Kesten's tree.

First, we  consider a functional  $A\,:\, \T_0 \longrightarrow  \N$ such
that $\{\bt;\,  A(\bt)\geq n\}$  is non  empty for  all $n\in
\N^*$. In the following theorems, we will add some assumptions on
$A$. These assumptions will vary from one theorem to another and in
fact we will consider three
different properties listed below (from the weaker to the stronger property):
for  all $\bt\in  \T_0$ and  all  leaf $x\in  \cl_0(\bt)$, there  exists
$n_0\in \N^*$  and $D(\bt,x)\geq 0$ (only  for the \ref{eq:+}  property) such
that  for all  $\bt'\in \T_0  $ satisfying  $A(\bt\circledast_x\bt')\geq
n_0$:
\begin{align}
\label{eq:>}   \tag{\text{Monotonicity}} 
 A(\bt\circledast_x\bt')& \geq A(\bt');\\
\label{eq:+}   \tag{\text{Additivity}} 
 A(\bt\circledast_x\bt')& =A(\bt')+D(\bt,x);\\
\label{eq:=}   \tag{\text{Identity}} 
 A(\bt\circledast_x\bt')& =A(\bt').
\end{align}
The  \ref{eq:=} property is  a particular case of the \ref{eq:+}
property with
$D(\bt,x)=0$;  and  the \ref{eq:+} property is a particular  case of the
\ref{eq:>} property.  We give examples of such functionals:
\begin{itemize}
\item  The  maximal  degree  $ M(\bt)=\max\{k_u(\bt),  \,  u\in  \bt\}$  has
  the  \ref{eq:=} property with $n_0=M(\bt)+1$.
\item The   cardinal   $  \sharp \bt   =\Card   (\bt)$ has
  the  \ref{eq:>} property with $n_0=0$ and has also the  \ref{eq:+}
  property with
  $n_0=0$ and $D(\bt,x)=\sharp \bt-1\geq 0$.  
   \item The  height $ H(\bt)=\max\{H(u),\, u\in \bt\}$
has the  \ref{eq:+} property with  $n_0=H(\bt)$ and $D(\bt,x)=H(x)\geq
0$.
 \end{itemize} 
Notice that the functional $A(\bt)=H(\bt)+ \sharp \bt$ does not
satisfy the \ref{eq:>}, \ref{eq:+} or \ref{eq:=} properties
as $A(\bt
\circledast_x\bt')= A(\bt') + A(\bt) - 1 - H(\bt')$ on $H(\bt') \leq
H(\bt ) -H(x)$.

\medskip



 We will   condition BGW trees with respect  to events $\A_n$
 of the form $\A_n=\{A(\tau)\ge n\}$ or $\A_n=\{A(\tau)=n\}$ or in order
 to avoid  periodicity arguments  $\A_n=\{A(\tau)\in [n,  n+n_1)\}$, for
 large  $n$. (Notice that the first two   cases boil  down to  the last one with
 respectively $n_1=+\infty $ and $n_1=1$.)

The next theorem states a general result concerning the local
convergence of {\bf critical} BGW tree conditioned on $\A_n$ toward
the Kesten's tree.
The proof of  this theorem is at  the end of
 this section.

 \begin{thrm}[Local convergence of critical BGW to Kesten's tree]
   \label{theo:main1}
  Let $\tau$ be  a {\bf critical} BGW tree with offspring  distribution
  $p$ such that $p(1)<1$
and
  let $\tau^*$ be a Kesten's tree associated with $p$. Let $\tau_n$ be a random
  tree    distributed    according    to   $\tau$    conditionally    on
  $\A_n=\{A(\tau)\in [n, n+n_1)\}$, where $n_1$ is fixed and   $\P(\A_n)>0$ for $n$ large enough. If one of the following conditions
  is satisfied:
\begin{enumerate}[(i)]
\item\label{it:>} $n_1=+\infty $ and $A$ satisfies the  \ref{eq:>} property;
\item\label{it:=}  $n_1\in
 \N^*\bigcup \{+\infty  \}$ and  $A$ satisfies the   \ref{eq:=} property;
\item\label{it:+}   $n_1\in \N^*\bigcup  \{+\infty  \}$,   $A$
  satisfies  \ref{eq:+} property and:
\begin{equation}
   \label{eq:rationA-n}
  \limsup_{n\to+\infty}\frac{\P(\A_{n+1})}{\P(\A_n)}\leq 1,
\end{equation}
\end{enumerate}
then, we have: 
\[
\tau_n\; \xrightarrow[n\rightarrow \infty ]{\textbf{(d)}}   \tau^*. 
\]
\end{thrm}


\begin{rmrk}[On the condition~\eqref{eq:rationA-n}]
   \label{rem:lim-inf}
   Assume    that   the    \ref{eq:+}    property    holds   and    that
   $\tau_n\;  \xrightarrow[n\rightarrow \infty  ]{\textbf{(d)}} \tau^*$.
   As $\T(\bt,  x)$ is open and  closed, we deduce from  the Portmanteau
   theorem                                                          that
   $\lim_{n\rightarrow  \infty }  \P(\tau_n\in \T(\bt,  x))=\P(\tau^*\in
   \T(\bt, x))$ for all $\bt\in  \T_0$ and $x\in \cl_0(\bt)$.  Thanks to
   \eqref{eq:tmp} below  in the  critical case  (with $m=1$),  we deduce
   that  $\lim_{n\rightarrow  \infty}  \P(\A_{n-D(\bt,  x)})/\P(\A_n)=1$
   provided  that  $\P(\tau^*\in  \T(\bt,   x))>0$ or, equivalently, that $\P(\tau\in  \T(\bt,   x))>0$ by Equation \eqref{eq:tau*-in-T}.   
   
   We  say  that  the
   functional  $D$ has  period  $d$ (with respect to $p$) if  the  smallest  group  in  $\Z$
   containing:
   \[
     \{D(\bt,     x)\,     \colon  \,    \bt\in     \T_0     ,\,     x\in
     \cl_0(\bt)\,\text{ and }\,   \P(\tau\in  \T(\bt,   x))>0\}
   \]
   is
   $d\Z$,  and that  $D$ is  aperiodic if  $d=1$.  Assuming  either that
   \eqref{eq:rationA-n} holds  or that the functional  $D$ is aperiodic,
   we deduce that:
\[
  \lim_{n\to+\infty}\frac{\P(\A_{n+1})}{\P(\A_n)}=1.
\]
Notice, one can also change the parameter $n_1$ so that the above convergence holds even if $d>1$.  

\end{rmrk}

\begin{rmrk}[Tail conditioning]
   \label{rem:tail}
   Let $\tau$ be a critical BGW tree with offspring distribution $p$ and
   let $\tau^*$ be a Kesten's tree associated with $p$.  For simplicity,
   let  us assume  that  $\P(A(\tau)=n)>0$ for  $n$  large enough.   Let
   $\tau_n$   be  a   random  tree   distributed  according   to  $\tau$
   conditionally on $\{A(\tau)=n\}$ and assume that:
\[
\tau_n\; \xrightarrow[n\rightarrow \infty ]{\textbf{(d)}}   \tau^*. 
\]
Since the distribution of $\tau$ conditionally on $\{A(\tau)\geq n\}$ is
a  mixture of the distributions of $\tau$ conditionally of
$\{A(\tau)=k\}$ for $k\geq n$, we deduce that $\tau$ conditionally on 
$\{A(\tau)\geq n\}$ converges in distribution toward $\tau^*$. 
In particular, as far as Theorem \ref{theo:main1} is concerned,  the
cases $n_1$ finite are the most delicate cases. 
\end{rmrk}

\begin{rmrk}[On $\P(\A_n)>0$]
  \label{rem:period}
Recall that by assumption,  $\{\bt;\,  A(\bt)\geq n\}$  is non  empty for  all $n\in
\N^*$. This is however not enough to ensure that $\P( A(\tau)\geq n)>0$ for $n$
large enough.
\end{rmrk}

 There is an extension  of  \ref{it:+}
 for a very special case  in the sub-critical case. 

 \begin{thrm}[Local convergence of sub-critical BGW to Kesten's tree]
   \label{theo:main2}
  Let $\tau$ be  a sub-critical BGW tree with  offspring distribution $p$
  with positive mean $m\in (0, 1)$,  and let $\tau^*$ be a  Kesten's
  tree associated with  $p$.  Let $\tau_n$ be a  random tree distributed
  according to $\tau$ conditionally  on $\A_n=\{A(\tau)\in [n, n+n_1) \}$
  with  $n_1\in \N^*  \cup \{+\infty  \}$  fixed, where  we assume  that
  $\P(\A_n)>0$ for $n$  large enough. If $A$  satisfies the \ref{eq:+}
  property  with
  $D(\bt,x)=H(x)$ and:
\begin{equation}
   \label{eq:rationA-n-2}
  \limsup_{n\to+\infty}\frac{\P(\A_{n+1})}{\P(\A_n)}\leq m,
\end{equation}
then, we have: 
\[
\tau_n\; \xrightarrow[n\rightarrow \infty ]{\textbf{(d)}}   \tau^*.
\]
\end{thrm}

\begin{rmrk}
The condition $D(\bt, x)=H(x) $ is very restrictive and holds essentially for
$A(\bt)=H(\bt)$ as we will see in the next section.
\end{rmrk}

\begin{proof}[Proof of Theorems \ref{theo:main1} and \ref{theo:main2}]
  As  we only consider  critical or subcritical  trees, the
  trees  $\tau_n$  belong  a.s.\  to  $\T_0$.  Moreover,  by  definition,
  a.s.\ the Kesten's tree belongs to $\T_1$. Therefore we can use Proposition
  \ref{prop:conv_dist} to  prove the convergence in  distribution of the
  theorems.

Let $n_1\in \N^* \cup \{ +\infty \}$ and set $\A_n=\{A(\tau) \in [n, n+n_1)\}$ in
  order to cover all the different cases of the two theorems. 
Let $\bt\in\T_0$. We have:
\[
\P(\tau_n=\bt)=\frac{\P(\tau=\bt,\ \A_n)}{\P(\A_n)}\le \frac{1}{\P(\A_n)}\ind_{\{A(\bt)\in[n,n+n_1)\}}.
\]
As $A(\bt)$ is
finite since $\bt\in \T_0$, we have $\ind_{\{A(\bt)\in[n,n+n_1)\}}=0$
  for $n> A(\bt)$. We deduce that:
\begin{equation}
   \label{eq:tn=t}
\lim_{n\rightarrow+\infty } \P(\tau_n=\bt)= 0=\P(\tau^*=\bt),
\end{equation}
as $\tau^*$ is a.s.\ infinite. This gives condition~\ref{it:cv-T0} of Proposition
\ref{prop:conv_dist}.

To obtain  the   convergence   in
distribution of  the sequence  $(\tau_n, n\in  \N^*)$ to  $\tau^*$, it
is enough to check  condition~\ref{it:cv-Ttx}  of Proposition
\ref{prop:conv_dist}, that is, for  $\bt\in\T_0$ and 
$x\in\cl_0(\bt)$ a leaf of $\bt$:
\begin{equation}
   \label{eq:liminfthm}
\liminf_{n\to+\infty}\P\bigl(\tau_n\in\T(\bt,x)\bigr)\geq
\P\bigl(\tau^*\in\T(\bt,x)\bigr). 
\end{equation}

\medskip

Let $\bt\in\T_0$ and  $x\in\cl_0(\bt)$ be a fixed leaf of $\bt$.  Since $\tau$ is
a.s.\ finite, we have:
\begin{equation}
   \label{eq:PA=PT}
\P(\tau\in\T(\bt,x),\A_n)= \sum_{\bt'\in\T_0} \P(\tau=\bt\circledast_x\bt')\ind_{\{n\le
  A(\bt\circledast_x\bt')<n+n_1\}}.
\end{equation}
Using \eqref{eq:Ptt'=PtPt'-v0}, \reff{eq:tau-in-T} and \reff{eq:tau*-in-T}, we get that for every
tree $\bt' \in \T$:
\begin{equation}
   \label{eq:Ptt'=PtPt'}
\P(\tau=\bt\circledast_x\bt')
=m^{H(x)} \P(\tau^*\in \T(\bt, x))\P(\tau=\bt').
\end{equation}
We deduce that:
\begin{equation}
   \label{eq:tn=t*}
\P(\tau\in\T(\bt,x),\A_n) 
=m^{H(x)} \P(\tau^*\in \T(\bt, x))\sum_{\bt'\in\T_0}\P(\tau=\bt')
\ind_{\{n\le
  A(\bt\circledast_x\bt')<n+n_1\}}.
\end{equation}
We now distinguish three different cases.


\medskip

\textbf{Assume $p$ critical, $n_1=+\infty$ and $A$ statisfies property \ref{eq:>}}.
 We deduce
from~\eqref{eq:PA=PT} that for $n\geq n_0$:
\[
\P(\tau\in\T(\bt,x),\A_n)=\sum_{\bt'\in\T_0}\P(\tau=\bt\circledast_x\bt')\ind_{\{n\le
  A(\bt\circledast_x\bt')\}}
\geq  \sum_{\bt'\in\T_0}\P(\tau=\bt\circledast_x\bt')\ind_{\{n\le
  A(\bt')\}},
\]
where we used the \ref{eq:>} property for the inequality.  Then,
using~\eqref{eq:Ptt'=PtPt'} with $m=1$, we get  that:
\[
\P(\tau\in\T(\bt,x),\A_n)
\geq  \P(\tau^*\in \T(\bt, x))\,\P\bigl(\A_{n}\bigr),
\]
which     gives    \reff{eq:liminfthm}.      This        proves         Theorem
\ref{theo:main1}~\ref{it:>}.

\medskip

\textbf{Assume $p$ critical and $A$ satisfies property \ref{eq:=}}.  In that case, we
have for $n\ge n_0$:
\[
\sum_{\bt'\in\T_0}\P(\tau=\bt')
\ind_{\{n\le
  A(\bt\circledast_x\bt')<n+n_1\}}=\sum_{\bt'\in\T_0}\P(\tau=\bt')
\ind_{\{n\le
  A(\bt')<n+n_1\}}=\P(\A_n),
\]
and we
obtain from \reff{eq:tn=t*}, with $m=1$, that
$
\P\bigl(\tau_n\in\T(\bt,x)\bigr)=\P\bigl(\tau^*\in\T(\bt,x)\bigr) 
$ for $n\geq n_0$. This yields~\eqref{eq:liminfthm} and 
 proves Theorem~\ref{theo:main1}~\ref{it:=}. 

\medskip

\textbf{Assume $A$ satisfies property \ref{eq:+}}. We deduce from \reff{eq:tn=t*}
that for $n\geq n_0$,  
\begin{align*}
\P(\tau \in\T(\bt,x),\A_n)
 &=m^{H(x)} \P(\tau^*\in \T(\bt, x))\sum_{\bt'\in\T_0}\P(\tau=\bt')
\ind_{\{n\le A(\bt')+D(\bt,x)<n+n_1\}}\\
& =m^{H(x)} \P(\tau^*\in \T(\bt, x))\, \P\bigl(n-D(\bt,x)\le
A(\tau)<n-D(\bt,x)+n_1\bigr)\\
&=m^{H(x)} \P(\tau^*\in \T(\bt, x))\,\P\bigl(\A_{n-D(\bt,x)}\bigr).
\end{align*}
Finally, we get: 
\begin{equation}\label{eq:tmp}
  \P\bigl(\tau_n\in\T(\bt,x)\bigr)=\frac{\P(\tau\in\T(\bt,x),\A_n)}{\P(\A_n)}
=m^{H(x)}\P\bigl(\tau^*\in\T(\bt,x)\bigr)
\frac{\P\bigl(\A_{n-D(\bt,x)}\bigr)}{\P(\A_n)}\cdot  
\end{equation}
If   \reff{eq:rationA-n}  holds   in  the   critical  case   ($m=1$)  or
\reff{eq:rationA-n-2} and  $D(\bt, x)=H(x)$ in the  sub-critical case, we
obtain~\eqref{eq:liminfthm}.          This        proves         Theorem
\ref{theo:main1}~\ref{it:+} and Theorem \ref{theo:main2}.
\end{proof}

\subsection{Applications to various conditioning}
\label{sec:appli}
We assume that the offspring distribution $p$
satisfies $p(1)<1$ and is critical ($m=1$) or sub-critical ($m<1$).

\subsubsection{\textbf{The height (Kesten's theorem)}}
\label{sec:Kesten}

We give a proof of Kesten's theorem, see Theorem
\ref{thm:cv-k1}. We consider the  functional of the tree given by its
height:   $A(\bt)=H(\bt)$. It satisfies the \ref{eq:+}  property, with
$n_0=H(\bt)+1$, as  for every  tree
$\bt\in\T_0$, every leaf $x\in\cl_0(\bt)$  and every
$\bt'\in \T_0$ such that $H(\bt\circledast _x\bt')\geq H(\bt)+1$, we have:
\[
H(\bt\circledast _x\bt')=H(\bt')+ H(x).
\]

We give a preliminary result. 

\begin{lmm}
   \label{lem:height}
   Let  $\tau$ be  a critical  or  sub-critical BGW  tree with  offspring
   distribution $p$ satisfying $p(1)<1$ with mean $m\leq 1$. Let
   $n_1\in \N ^*  \cup \{ +\infty \}$.  Set $\A_n=\{A(\tau)\in  [n, n+n_1)\}$ for
   $n\in \N^*$. We have:
\[
\lim_{n\to+\infty}\frac{\P(\A_{n+1})}{\P(\A_n)}=m.
\]
\end{lmm}

We then deduce the following corollary as a direct
consequence of  Theorem~\ref{theo:main1}~\ref{it:+} in the critical case or of Theorem
\ref{theo:main2} in the sub-critical case.

\begin{crllr}[Kesten's theorem]
   \label{cor:kesten-H}
Let $\tau$ be  a  critical or sub-critical  BGW tree with offspring
distribution $p$ satisfying $p(1)<1$, and
  $\tau^*$ a Kesten's tree associated with $p$. 
Let $\tau_n$ be a random
  tree    distributed    according    to   $\tau$    conditionally    on
  $\{H(\tau)=n\}$ (resp. $\{H(\tau)\geq n\}$). 
Then, we have: 
\[
\tau_n\; \xrightarrow[n\rightarrow \infty ]{\textbf{(d)}}   \tau^*. 
\]
\end{crllr}

\begin{proof}[Proof of Lemma \ref{lem:height}]
 We shall consider the case $n_1=1$, the other  cases being  deduced
 using  Remark \ref{rem:tail}. So, we have $\A_n=\{H(\tau)=n\}$.  Recall that for any tree $\bt$, $z_n(\bt)=\Card\{u\in\bt,\ H(u)=n\}$
  denotes 
the size of the $n$-th generation of the tree $\bt$. We set
$Z_n=z_n(\tau)$, so that  $(Z_n, {n\ge 0})$ is  a BGW
process. Notice that:
\[
\A_{n}=\{Z_{n+1}=0\}\setminus \{Z_{n} =0\}.
\]
Recall that $g$ denotes the generating function of the offspring distribution
$p$. Let  $g_n$ be the generating function of  $Z_n$. In particular,
we have $g_1=g$. Using the branching property of the BGW tree, we have
that $Z_{n+1}$ is distributed as $\sum_{i=1}^{k_\root(\tau)} z_n(\tau_i)$,
where $(\tau_i, i\in \N^*)$ are independent BGW tree with offspring
distribution $p$ and independent of $Z_1=k_\root(\tau)$. This gives:
\[
g_{n+1}(s)=\E\left[\prod_{i=1}^{Z_1}
  s^{z_n(\tau_i)}\right]
=\E\left[g_n(s) ^{Z_1}\right]=g(g_n(s)).
\]
We have $\P(\A_{n})=\P(Z_{n+1}=0) - \P(Z_{n}=0)=g_{n+1}(0) - g_{n}(0)$. 
Since $\tau$ is critical or sub-critical, it is a.s.\ finite and we deduce that
$\lim_{n\rightarrow+\infty } g_n(0)=\lim_{n\rightarrow+\infty }
\P(z_n(\tau)=0)=1$. We have:
\[
\frac{\P(\A_{n+1})}{\P(\A_n)}
=\frac{g_{n+2}(0) - g_{n+1}(0)}
{g_{n+1}(0)- g_n(0)}\cdot
\]
Using Taylor formula at $g_n(0)$, we get:
\begin{align*}
g_{n+2}(0) &
=g\big(g_n(0)+(g_{n+1}(0) - g_n(0))\big)\\
& =g_{n+1}(0) + (g_{n+1}(0) - g_n(0)) \,g'(g_n(0)) + o(g_{n+1}(0) - g_n(0)).
\end{align*}
This gives that:
\[
\frac{\P(\A_{n+1})}{\P(\A_n)}
=g'(g_n(0))+ o(1)\xrightarrow[n\rightarrow \infty ]{ } m.
\]
\end{proof}

\subsubsection{\textbf{The maximal out-degree, critical case}}
\label{sec:out}

Following \cite{h:cgwtmod}, we consider the functional of the tree given
by its maximal out-degree: $A(\bt)=M(\bt)$, with:
\[
M(\bt)=\sup_{u\in\bt}k_u(\bt).
\]
Notice it satisfies the  \reff{eq:=} property, with $n_0= M(\bt)+1$, 
 as  for every  tree
$\bt\in\T_0$, every leaf $x\in\cl_0(\bt)$ and every
$\bt'\in \T_0$ such that $M(\bt\circledast _x\bt')\geq M(\bt)+1$, we have:
\[
M(\bt\circledast _x\bt')=M(\bt').
\]

The next corollary is then a consequence of  Theorem
\ref{theo:main1}~\ref{it:=} with $n_1\in \{1, +\infty
\}$. For $n_1=1$, in the proof of this theorem, the condition
$\P(\A_n)>0$ for $n$ large enough is easily  replaced by the
convergence along the  sub-sequence $\{n; p(n)>0\}$ given by the support
of $p$. 

\begin{crllr}[Conditioning on the maximal out-degree]
   \label{cor:kesten-M}
   Let  $\tau$ be  a critical  BGW tree  with offspring  distribution $p$
   satisfying $p(1)<1$ and unbounded support, and
   let $\tau^*$ be a Kesten's  tree associated with $p$. Let  $\tau_n$ be a
   random  tree   distributed  according  to  $\tau$   conditionally  on
   $\{M(\tau)= n\}$ (resp.  $\{M(\tau)\geq n\}$).  Then, we have  along the sub-sequence
   $\{n; p(n)>0\}$ (resp. with $n\in \N^*$):
\[
\tau_n\; \xrightarrow[n\rightarrow \infty ]{\textbf{(d)}} \tau^*.
\]
\end{crllr}

Intuitively, the height of the largest node in the BGW tree $\tau_n$ goes to infinity, so it disappears in the local limit.

\subsubsection{\textbf{The Horton-Strahler number, critical case}}
\label{sec:horton}
The Horton-Strahler number of a finite tree was initially introduced in
hydrogeomorphism by Horton (1845) and redefined later by Strahler
(1952), see~\cite{viennot,esparza,dp:rftt,bdr:HS21}
and references therein. 
It is also known as  the register function
in computer science. This number measures the unbalance of a tree, and
was initially introduced for binary trees ($k_u(\bt)\in \{ 0, 2\}$
for $u\in \bt$). Recall the  definition of the  sub-tree $\cs_u(\bt)$ of $\bt$  above the
node  $u\in \bt$ given  in  \reff{eq:defSu}.   The Horton-Strahler number $\Sigma(\bt)$ of a finite binary
tree  $\bt$  is defined  recursively  by  $\Sigma(\{\root\})=0$ and,  if
$k_\root(\bt)=2$,                                                   then
$\Sigma(\bt)=\Sigma_1\vee  \Sigma_2   +  \ind_{\{\Sigma_1=\Sigma_2\}}$
where $\Sigma_u= \Sigma(\cs_u(\bt))$ for $u\in \{1, 2\}$.  \medskip

Its generalization, which we still write $\Sigma$,  to finite
general rooted tree also called the register function is defined
as follows, for $\bt \in \T_0$: 
\[
\Sigma(\bt)=
\begin{cases}
   0 & \text{ if } \bt=\{\root\},\\
\max\{\Sigma_{(1)}, \Sigma_{(2)}+1, \ldots, \Sigma_{(k)} +
k-1\} & \text{ if } k_\root(\bt)=k\geq 1,
\end{cases}
\]
where $\Sigma_{(1)}\geq  \cdots\geq  \Sigma_{(k)}$ is the
non-increasing reordering of $(\Sigma_1, \ldots, \Sigma_k)$ and
$\Sigma_u= \Sigma(\cs_u(\bt))$ for $u\in \{1, \ldots, k\}$. 
Clearly, we have that $\Sigma(\bt)\leq  \Sigma_{(1)} + k_\root(\bt) -1$. 
Notice also that the two definitions of $\Sigma$ coincide on binary trees.

\medskip

The  functional of the tree $A=\Sigma$
has the  \ref{eq:=} property, with $n_0=\Sigma(\bt)+M(\bt)-1$,
where $M(\bt)$ is the maximal out-degree of $\bt$. Indeed, for 
trees $\bt, \bt'\in\T_0$ and  leaf  $x\in\cl_0(\bt)$ of $\bt$ 
 such that $\Sigma(\bt\circledast _x\bt')\geq
 \Sigma(\bt)+M(\bt)-1$,
it is easy to check that:
\[
\Sigma(\bt\circledast _x\bt')=\Sigma(\bt').
\]
Notice that $\P(\Sigma(\tau)=n)>0$ for all $n\geq \inf\{k\geq 1\, \colon p(k+1)>0\}$. 
The next corollary is then a consequence of
Theorem~\ref{theo:main1}~\ref{it:=} with $n_1\in \{1, +\infty 
\}$.

\begin{crllr}[Conditioning on the Horton-Strahler number]
   \label{cor:kesten-Horton}
   Let  $\tau$ be  a critical  BGW tree  with offspring  distribution $p$
   satisfying  $p(1)<1$, 
    and  let   $\tau^*$  be a  Kesten's  tree
   associated  with  $p$. Let  $\tau_n$  be  a random  tree  distributed
   according  to $\tau$  conditionally on  $\{\Sigma(\tau)= n\}$  (resp.
   $\{\Sigma(\tau)\geq  n\}$).  Then,  we  have:
\[
\tau_n\; \xrightarrow[n\rightarrow \infty ]{\textbf{(d)}} \tau^*.
\]
\end{crllr}

\begin{rmrk}[Extension to CRT]
  \label{rem:sl-horton}
  A very recent result on scaling limit of the  Horton-Strahler number 
for a   BGW conditioned to be large appears in \cite{khanfir2024}. 
This presentation explains in particular why the fluctuations of the
Horton-Strahler number are asymptotically coupled with 
deterministic oscillations. 
\end{rmrk}

\subsubsection{\textbf{The largest generation,  critical case}}
\label{sec:large-gen}

We consider the functional of the tree given
by its largest generation: $A(\bt)=\cz(\bt)$ with:
\[
\cz (\bt)=\sup_{k\geq 0} z_k(\bt).
\]
Notice the functional has the  \ref{eq:>} property  as $\cz
(\bt\circledast _x\bt')\geq \cz (\bt')$. 
The next corollary is then a direct consequence of
Theorem~\ref{theo:main1}~\ref{it:>}; this result appears in \cite{xin-loc-cv}.

\begin{crllr}[Conditioning on the largest generation]
   \label{cor:kesten-Z}
Let $\tau$ be  a  critical   BGW tree with offspring
distribution $p$ satisfying $p(1)<1$,  and
  let $\tau^*$ be a Kesten's tree associated with $p$. 
Let $\tau_n$ be a random
  tree    distributed    according    to   $\tau$    conditionally    on
  $\left\{\cz (\tau)\geq n\right\} $. 
Then, we have:
\[
\tau_n\; \xrightarrow[n\rightarrow \infty ]{\textbf{(d)}} \tau^*.
\]
\end{crllr}
\begin{rmrk}[Conditioning on $\left\{\cz (\tau) =  n\right\} $]
   \label{rem:Z}
   Notice  that the  functional $\cz$  does not  satisfy the  \ref{eq:+}
   property.      Thus,    considering     the    conditioning     event
   $\left\{\cz (\tau) =  n\right\} $ is still an  open problem. However,
   as  noted  in \cite{xin-loc-cv},  if  the  support of  the  offspring
   distribution  is bounded,  then  the functional  $\cz$ satisfies  the
   \ref{eq:=} property  (take $n_0=\cz(\bt)  + \alpha^{H(\bt)  - H(x)}+1$
   with   $\alpha=    \sup   \supp(p)$    being   finite),    and   then
   Corollary~\ref{cor:kesten-Z}      holds      with      the      event
   $\left\{\cz      (\tau)\geq     n\right\}      $     replaced      by
   $\left\{\cz  (\tau)=  n\right\}  $  and  considering  the  (infinite)
   subsequence of $n$ for which this event has positive probability.
\end{rmrk}

\subsubsection{\textbf{The total progeny, critical case}}
\label{sec:total}

The local convergence  in distribution of  the critical tree,  conditionally on
its  total size  being  large,  to the  Kesten's  tree  appears implicitly  in
\cite{k:gwctp} and was first explicitly stated in \cite{ap:tvmcdgwp}. We
give   here  an   alternative  proof.    We  consider   the  functional:
$A(\bt)=\sharp \bt$, with  $\sharp \bt=\Card(\bt)$ the total size  of $\bt$, which
has the   \ref{eq:+} property  as
  for    every   trees    $\bt,\bt'\in\T_0$   and    every   leaf
$x\in\cl_0(\bt)$:
\[
\sharp \bt\circledast _x\bt'=\sharp \bt'+\sharp \bt-1.
\]
 
Recall $d$ is the period of  the
offspring distribution $p$.  Using B\'ezout's identity, it is  classical  to
  check that:
  \begin{equation}
   \label{eq:|t|-period}
   \P(\sharp \tau\in 1+ d\N)=1
     \end{equation}
     and that  there exists $n_0\in \N^*$  such that $\P(\sharp \tau=1+dn)>0$
     for all $n\geq n_0$.
The next  lemma is a direct  consequence of
Dwass formula  and the strong ratio  limit theorem. Its proof  is at the
end of this section.

\begin{lmm}
   \label{lem:total-size}
Let $\tau$ be  a  critical   BGW tree with offspring
distribution $p$ satisfying $p(1)<1$ with period $d$. Let $n_1\in \{d, +\infty
\}$ and set $\A_n=\{\sharp \tau \in [n, n+n_1)\}$. Then, we have:
\[
\lim_{n\to+\infty}\frac{\P(\A_{n+1})}{\P(\A_n)}=1.
\] 
\end{lmm}

We then deduce the following corollary as a direct
consequence of  Theorem~\ref{theo:main1}~\ref{it:+}.

\begin{crllr}[Conditioning on the total size]
   \label{cor:kesten-total}
Let $\tau$ be  a  critical   BGW tree with offspring
distribution $p$ satisfying $p(1)<1$ with period $d$, and
  let $\tau^*$ be a Kesten's tree associated with $p$. 
Let $\tau_n$ be a random
  tree    distributed    according    to   $\tau$    conditionally    on
  $\{\sharp \tau=1+dn\}$ for $n$ large enough (resp. on  $\{\sharp \tau\geq n\}$). 
Then, we have: 
\[
\tau_n\; \xrightarrow[n\rightarrow \infty ]{\textbf{(d)}}   \tau^*. 
\]
\end{crllr}

As in Section~\ref{sec:out}, we may observe in $\tau_n$ a large node whose height and size go to infinity and which disappears in the local limit, see~\cite{kr:cond-cauchy} in this direction.


The   end  of   the  section   is  devoted   to  the   proof  of   Lemma
\ref{lem:total-size}.    
We first recall  Dwass formula from \cite{d:tpbprrw} that links the distribution of
the total progeny of BGW trees to the distribution of 
random walks and then the strong ratio limit theorem. 

  Let $(\zeta_k,  k\in \N^*)$  be a
  sequence  of independent  random  variables  distributed according  to
  $p$. Set $S_n=\sum_{k=1}^n  \zeta_k$ for $n\in \N^*$.



\begin{lmm}[Dwass formula]\label{lem:Dwass}
  Let  $\tau$ be  a BGW  tree  with   offspring
  distribution $p$. 
  Then, for every  $n\ge 1$, we have:
  \[
\P\left(\sharp \tau=n\right)
=\frac{1}{n}\P(S_n =n-1).
\]
\end{lmm}



We also recall the strong ratio limit theorem that can be found for
instance in  \cite{n:tece},  see also \cite{s:prw} Theorem T1 p.49.

\begin{lmm}[Strong ratio limit theorem]\label{lem:strong_ratio}
 Assume that $p$ is  critical aperiodic with $p(1)<1$. Then, we  have, for every
  $\ell\in\Z$,
\[
\lim_{n\to+\infty}\frac{\P(S_n=n+\ell)}{\P(S_n=n)}=1\quad\text{and}\quad
\lim_{n\to+\infty}\frac{\P(S_{n+1}=n+1)}{\P(S_n=n)}=1.
\]
\end{lmm}

  \begin{rmrk}[Strong ratio limit theorem in the periodic case]
  We first stress that the strong ratio theorem in \cite{s:prw,n:tece} requires
  $p$ to  be strongly aperiodic,  that is, the smallest  sub-group $G_k$
  which contains $k+\supp(p)$ is $\Z$  for all $k\in \Z$. Taking $k=0$,
  we get  that $p$  is aperiodic  if it is  strongly aperiodic.   On the
  other hand,  if $p$ is aperiodic  and $1>p(0)>0$, then we  deduce that
  $k$ belongs  to $G_k$, and thus  $-k$ also; this readily  implies that
  $\supp(p)\subset G_k$ and  thus $G_k=\Z$, giving that  $p$ is
  strongly aperiodic.

 Eventually, if $p$ has period $d$, we refer to Equation~(4.10) in 
 \cite{ad:llcgwtisc} for the extension of the strong ratio limit theorem. 
  \end{rmrk}

\begin{proof}[Proof of Lemma \ref{lem:total-size}]
  We  shall  assume  for  simplicity  that $p$  is  aperiodic  and  take
  $n_1=1$. The cases  $d\geq 2$ or $n_1=+\infty $  are left to the reader.
  Using Dwass formula, we have:
\[
\frac{\P(\sharp \tau=n+1)}{\P(\sharp \tau=n)}=\frac{n}{n+1}\frac{\P(S_{n+1}
  =n)}{\P( S_n =n-1)}\cdot
\]
Using  the strong ration limit theorem, we get:
\[
\lim_{n\to+\infty}\frac{\P(\sharp \tau=n+1)}{\P(\sharp \tau=n)} 
=\lim_{n\to+\infty}\frac{\P(S_{n+1}=n)}{\P(S_n=n-1)}
=1.
\]
This ends the proof of   Lemma
\ref{lem:total-size}.
\end{proof}

\subsubsection{\textbf{The number of leaves, critical case}}
\label{sec:L0}

Recall $\cl_0(\bt)$  denotes the set of  leaves of a tree  $\bt$ and set
$L_0(\bt)$  for its  cardinal. We  shall  consider a  critical BGW  tree
$\tau$   conditioned   on   the   number  of   leaves,   that   is,   on
$\{L_0(\tau)=n\}$.    Such    a    conditioning   appears    first    in
\cite{ck:rncpccgwta} with a  second moment condition. We  prove here the
convergence in distribution of the  conditioned tree to the  Kesten's tree in
the  critical  case  without  any additional  assumption  using  Theorem
\ref{theo:main1}.

The functional $A(\bt)=L_0(\bt)$ satisfies  the \ref {eq:+} property
with $n_0=1$ and $D(\bt, x)=L_0(\bt) -1$ as 
for every trees $\bt,\bt'\in\T_0$ and every leaf $x\in\cl_0(\bt)$:
\[
L_0(\bt\circledast _x\bt')=L_0(\bt')+ L_0(\bt)-1.
\]

The  next  lemma  due  to Minami  \cite{m:nvgdgwt}  gives  a  bijection 
 between the  leaves of a finite tree $\bt$  and the nodes
of a tree $\bt_ {\{0\}}$. Its proof is given at the end of this section.

\begin{lmm}[The tree of leaves]
   \label{lem:t0}
   Let  $\tau$ be  a critical  BGW tree  with offspring  distribution $p$
   such that $p(1)<1$.  Then   $L_0(\tau)$  is  distributed  as
   $\sharp \tau_ {\{0\}}$, where  $\tau_ {\{0\}}$ is a critical  BGW tree
   with 
   offspring distribution  $p_{(0)}$ such that $p_{(0)}(1) <1$.
\end{lmm}


It is elementary to check from the construction of $\tau_{\{0\}}$ in the
proof of Lemma~\ref{lem:t0} below
that the period $d_ {0}$
of $p_ {\{0\}}$ is given by:
\[
  d_0=\max\{ k;\,  \supp(p)\subset k \N +1\}.
\]
This,  Lemma~\ref{lem:t0} and~\eqref{eq:|t|-period}  imply in  turn that
$\P(L_{0}(\tau)\in   1+  d_0   \N)=1$   and   that  there   exists
$n_0\in   \N^*$  such   that   $\P(L_0(\tau)=1+d_0 n)>0$  for   all
$n\geq n_0$.   The following corollary  is then a direct  consequence of
Lemma \ref{lem:total-size} and Theorem \ref{theo:main1}~\ref{it:+}.

\begin{crllr}[Conditioning on the number of leaves]
   \label{cor:kesten-leaves}
   Let $\tau$  be a  critical BGW tree  with offspring  distribution $p$,
  such that $p(1)<1$,  and $\tau^*$  a Kesten's  tree associated
   with  $p$.     Let
   $\tau_n$   be  a   random  tree   distributed  according   to  $\tau$
   conditionally  on  $\left\{L_0(\tau)=1+ d_0 n\right\}$
   (resp.  $\{L_0(\tau)\geq n\}$).  Then, we have:
\[
\tau_n\; \xrightarrow[n\rightarrow \infty ]{\textbf{(d)}}   \tau^*. 
\]
\end{crllr}


\begin{proof}[Proof of Lemma \ref{lem:t0}]
We first describe the correspondence  of Minami.  The left-most leaf (in
the  lexicographical order)  of $\bt$  is mapped  on the  root of  $\bt_
{\{0\}}$. In the example of  Fig.~\ref{fig:t0:root}, the leaves of the
tree $\bt$ are labeled from 1 to 9, and the left-most leaf is 1.


\begin{figure}[H]
  \begin{tikzpicture}
  [level 3/.style={sibling distance=10mm}, level 1/.style={sibling distance=20mm}]
  \node[inner sep=-.6,outer sep=-.6] {} [grow=north]
child {node[inner sep=-.6,outer sep=-.6]{}
	child{node[inner sep=-.6,outer sep=-.6] (c) {}
		child{node{9}}
		child{node  {8}}
		child{node{7}}
              }
              child{node (a) {6}}
	}
child {node (a) {5}}
child {node [inner sep=-.6,outer sep=-.6]{}
  child{node[inner sep=-.6,outer sep=-.6]  {}
    child {node{4}}
    child{node (b) {3}}
    child{node{2}}}
	child{node{1}}
      }
      ;
      
\end{tikzpicture}
\caption{A finite tree $\bt$.}
\label{fig:t0:root}
\end{figure}

Then consider all the subtrees that are attached to the branch between
the root and the left-most leaf. All the left-most leaves of these
subtrees are mapped on the children of the root of $\bt_ {\{0\}}$, they form
the population at generation 1 of the tree $\bt_ {\{0\}}$. In  Fig.~\ref{fig:t0:G1}, 
the considered sub-trees are surrounded by dashed lines, and the
leaves at generation 1 are labeled $\{2,5,6\}$. Remark that the
sub-tree that contains the leaf 5 is reduced to a single node (this
particular leaf).

\begin{figure}[H]
  \begin{tikzpicture}
  [level 3/.style={sibling distance=10mm}, level 1/.style={sibling distance=20mm}]
  \node[inner sep=-.6,outer sep=-.6] {} [grow=north]
child {node[inner sep=-.6,outer sep=-.6]{}
	child{node[inner sep=-.6,outer sep=-.6] (c) {}
		child{node{9}}
		child{node  {8}}
		child{node{7}}
              }
              child{node (a) {6}}
	}
child {node (a) {5}}
child {node [inner sep=-.6,outer sep=-.6]{}
  child{node[inner sep=-.6,outer sep=-.6]  {}
    child {node{4}}
    child{node (b) {3}}
    child{node{2}}}
	child{node{1}}
      }
      ;
      
      \draw [dashed] (a) ellipse (.4cm and .3cm);
      \draw [dashed] (b) ellipse  (1.4cm and 1.5cm);
       \draw [dashed] (c) ellipse  (2.6cm and 1.7cm);
     \end{tikzpicture}
\caption{The sub-trees attached to the branch between the root and the
  leaf labeled 1.}
\label{fig:t0:G1}
\end{figure}

Then perform the  same procedure inductively at each  of these sub-trees
to construct the tree $\bt_ {\{0\}}$. In  Fig.~\ref{fig:t0}, we give the tree
$\bt_ {\{0\}}$ associated with the tree $\bt$.

\begin{figure}[H]
  \begin{tikzpicture}
  \node {1} [grow=north]
child {node{6}
	child{node {7}
		child{node{9}}
		child{node  {8}}
              }
	}
child {node  {5}}
child {node {2}
  child{node  {4}}
	child{node{3}}
      }
      ;
     \end{tikzpicture}
     \caption{The tree $\bt_ {\{0\}}$ associated with the tree $\bt$ from Fig.~\ref{fig:t0:root}.}
\label{fig:t0}
\end{figure}

Using the branching property, we get that all the sub-trees that are
attached to the left-most branch of a BGW tree $\tau$ are independent and
distributed as $\tau$. Therefore, the tree $\tau_ {\{0\}}$ is still a
BGW tree. 

Next,  we   compute  the   offspring  distribution,  $p_   {\{0\}}$,  of
$\tau_ {\{0\}}$.  Since $p$ is  critical and $p(1)<1$, we  get $p(0)>0$.
Let  $N$ denote the generation of the left-most leaf. It is easy to see
that  this  random variable  is  distributed  according to  a  geometric
distribution with parameter $p(0)$, that is, for every $n\ge 1$:
\[
\P(N=n)=\bigl(1-p(0)\bigr)^{n-1}p(0).
\]
Let $\zeta$ be a random variable  with distribution $p$ and mean $m$,
and let $X$ be distributed as   $\zeta$
conditionally on $\{\zeta>0\}$, that is,
$\P(X=n)=p(n)/(1-p(0))$ for every $n\ge 1$:
In particular, we have:
\[
\E[X]=\frac{m}{1-p(0)}\cdot
\]

We
denote by $(X_1,\ldots,X_{N-1})$ the respective numbers of offsprings of
the nodes  on the  left-most branch (including  the root,  excluding the
leaf). Then,  using again  the branching  property, these  variables
are, conditionally on $N$, 
independent  and   distributed   as   $X$. 
Thus, the number of children of the root in the tree $\tau_ {\{0\}}$ is the
number of the sub-trees attached to the left-most branch that is:
\[
\zeta'=\sum_{k=1}^{N-1}(X_k-1).
\]
By construction $p_{\{0\}}$ is the probability distribution of
$\zeta'$. Since $p(0)>0$, we get that  that $\P(N=
1)> 0$ and thus  $p_{\{0\}}(0)>0$. 
We now compute the mean of  $p_{\{0\}}$:
\[
\E\left[\zeta'\right]  =\E[N-1]\E[X-1]
= \left(\frac{1}{p(0)}-1\right)\left(\frac{m}{1-p(0)}-1\right)
=\frac{1
}{p(0)}\Bigl(m-\bigl(1-p(0)\bigr)\Bigr).
\]
In particular, if the BGW tree $\tau$ is critical ($m=1$), then
$ \E\left[\zeta'\right]=1$, and thus the BGW tree  $\tau_{\{0\}}$ is also
critical.
\end{proof}

\subsubsection{\textbf{The  number of nodes with given out-degree,
  critical case}}
\label{sec:out-d}
 
 Let
$\ca$ be a non-empty subset of $\N$ and for a tree $\bt$, we define the subset of
nodes with out-degree in $\ca$:
\[
\cl_\ca(\bt)=\{u\in\bt,\ k_u(\bt)\in\ca\}
\]
and $L_\ca(\bt)=\Card\bigl(\cl_\ca(\bt)\bigr)$ its cardinal. 
(If $\ca=\N$, then $L_\ca(\bt)$ is the total number of nodes of $\bt$;
and if $\ca=\{0\}$, then $L_\ca(\bt)=L_0(\bt)$ is the total number of
leaves of $\bt$.)
The functional $L_\ca$  satisfies  Property
\reff{eq:+} with  $D(\bt, x)=L_\ca(\bt)
-\ind_{\{0\in\ca\}}\geq 0$, that is, for $\bt, \bt'\in \T_0$ and $x\in \cl_0(\bt)$:
\[
L_\ca(\bt\circledast_x\bt')=L_\ca(\bt')+L_\ca(\bt)-\ind_{\{0\in\ca\}}.
\]
Moreover, when $\bt\in T_0$ is such that $L_A(\bt)>0$, there  exists
a  bijection,  generalizing  Minami's correspondence,  between  the  set
$\cl_\ca(\bt)$ and   a  tree   $\bt_\ca$  so   that
$L_\ca(\bt)=\sharp  \bt_A$, see  Rizzolo \cite{r:slmbtgwtcnvodgs}.   Let
$\tau$  be  a  BGW  tree  with  offspring  distribution  $p$  such  that
$p(\ca)>0$, where we set:
\[
p(\ca)=\sum_{n\in\ca}p(n).
\]
Without   further   assumption   on   $p$,   notice   that   the   event
$\{L_\ca(\tau)>0\}$ has positive probability and that:
\begin{equation}
  \label{eq:bgw-in-T0}
  \P(\tau\not \in \T_0, \, L_\ca(\tau)<+\infty )=0.
\end{equation}
The latter equality is obvious if  $p$ is sub-critical or critical. When
$p$ is super-critical with
$p(0)>0$ (and thus $q>0$), it is an  immediate consequence of
Corollary~\ref{cor:TS} using  the strong  law of  large number  and that
there are  an infinity of  independent BGW sub-critical  sub-trees, with
offspring distribution  $\tilde p$  defined by  \reff{eq:tilde-p}, which
are  grafted  on  the  infinite  backbone of  the  individuals  of  type
survivor. On the other hand, if $p(0)=0$, then the number of children of the individuals in the infinite left branch of $\tau$, that is $(k_u(\tau), u\in \cup_{\ell\in \N} \{1\}^\ell)$ are independent random variables
distributed as $p$; thus a.s.\ $L_\ca(\tau)$ is infinite.



\begin{lmm}[Theorem 6 in \cite{r:slmbtgwtcnvodgs}]
  \label{lem:tA}
  Let $\tau$  be a BGW  tree with  offspring distribution $p$  such that
  $p(\ca)>0$, and $p$ is either  critical with $p(1)<1$ or sub-critical.
  Then, there exists a BGW  tree $\tau_ \ca$ with offspring distribution
  $p_    {\ca}$    such    that    $L_\ca(\tau)$,    conditionally    on
  $\{L_\ca(\tau)>0\}$,   is  distributed   as  $\sharp   \tau_  \ca   $.
  Furthermore  if  $p$ is critical then $p_\ca$ is critical with $p_\ca(1)<1$,
  and $p_\ca$ is sub-critical otherwise.
\end{lmm}


According to~\cite{r:slmbtgwtcnvodgs},
the period $d_\ca$ of  $p_\ca$
depends only on the sets $\ca$ and $\supp(p)$. 

\medskip

The next  result is  a generalization of  Section~\ref{sec:L0}. It  is a
direct   consequence   of   Lemma   \ref{lem:total-size}   and   Theorem
\ref{theo:main1}~\ref{it:+}.

\begin{crllr}[Conditioning on the number of nodes with given out-degree]
   \label{cor:kesten-ca}
   Let  $\tau$ be  a critical  BGW tree  with offspring  distribution
   $p$, such that  $p(\ca)>0$ and $p(0)>0$, and let $\tau^*$ be Kesten's  tree associated
   with $p$.  Let  $\tau_n$ be  a random
   tree    distributed   according    to    $\tau$   conditionally    on
   $\left\{L_\ca(\tau)= 1+ d_     \ca \, n \right\}$     (resp.
   $\{L_\ca(\tau)\geq n\}$).  Then, we have:
\[
\tau_n\; \xrightarrow[n\rightarrow \infty ]{\textbf{(d)}}   \tau^*. 
\]
\end{crllr}

\subsection[Number of nodes with given out-degree, sub-critical case]{Conditioning on the number of nodes with given out-degree, sub-critical case}\label{sec:sub-critical}

Theorem \ref{theo:main2} deals with sub-critical offspring distributions
and  applies  essentially   for  the  conditioning  on   the  height,  see
Corollary~\ref{cor:kesten-H}.   We  complete   the  picture  of  Theorem
\ref{theo:main1} studying the local limit of BGW trees by conditioned on
$\{L_\ca(\tau)=n\}$ in the sub-critical case.

We shall first present in Section~\ref{sec:def-ca} a family of offspring
distributions for which the conditioned trees have the same
distribution. If this family contains one critical 
offspring
distribution, which is the so-called \emph{generic} case, then we can
use the results on the critical case from Section~\ref{sec:out-d}. 
When this is not the case, this is the so called \emph{non-generic}
  case. Then the local limit exhibit a \emph{condensation} phenomenon,
  that is, there exists a (unique) node in the random limiting tree with
  infinite out-degree, see Proposition~\ref{prop:non-gene-sub}. This
  result is proved in~\cite{ad:llcgwtcc} and is more technical as the
  limit does not belong to $\T_0\cup \T_1$, and thus one needs a new
  characterization of the local convergence. 
Eventually, we discuss the generic/non-generic properties depending on
the set $\ca$ isn Section~\ref{sec:gene-non-gene}.

\subsubsection{An equivalent probability}
\label{sec:def-ca}
Let $p$ be an offspring distribution and $\ca\subset \N$. Recall
$p(\ca)=\sum_{n\in\ca}p(n)$. 
We assume that $p(\ca)>0$ and we define:
\[
  I_\ca=\Big\{\theta\in    (0,   +\infty    )\,   \colon\,    \sum_{k\in
    \ca}\theta^{k-1}p(k)<+\infty\quad                    \text{and}\quad
  \sum_{k\not\in\ca}\theta^{k-1}p(k)< 1\Big\}.
\]
For $\theta\in I_\ca$, we set for every $k\in\N$:
\[
p_\theta^\ca(k)=\begin{cases}
\theta ^{k-1}p(k) & \mbox{if }k\not\in\ca,\\
c_\ca(\theta)\, \theta^{k-1} p(k) & \mbox{if }k\in\ca,
\end{cases}
\]
where $c_\ca(\theta)$ is a finite positive constant that makes $p_\theta^\ca$ a
probability measure on $\N$, namely:
\[
c_\ca(\theta)=\frac{1-\sum_{k\not\in\ca}\theta^{k-1}p(k)}
{\sum_{k\in\ca}\theta^{k-1}p(k)}\cdot   
\]
Remark  that  $I_\ca$   is  exactly  the  set  of   $\theta$  for  which
$p_\theta^\ca$ is indeed a probability  measure  with the same support as
$p$.  Remark also  that  $I_\ca$  is an  interval  that  contains 1,  as
$p_1^\ca=p$.

The  following proposition gives the connection between $p$ and
$p^\ca_\theta$.  It  generalizes  the  results  already
  obtained for  the total progeny,  $\ca=\N$, in  \cite{k:gwctp}  and for
  the number  of leaves, $\ca=\{0\}$, in~\cite{adh:pgwttvmp}. 
Since $p$ and $p^\ca_\theta$ have the same support, we deduce that 
$\P\bigl(L_\ca(\tau)=n\bigr) >0$ if and only if
$\P\bigl(L_\ca(\tau_{[\theta]})=n\bigr) >0$.
  
\begin{prpstn}[Same conditional distribution]
\label{prop:prob_equiv}
Let  $\tau$  be  a  BGW   tree  with  offspring  distribution  $p$.  Let
$\ca\subset  \N$ such  that $p(\ca)>0$  and let  $\theta\in I_\ca$.  Let
$\tau_{[\theta]}$   be   a   BGW  tree   with   offspring   distribution
$p_\theta^\ca$.   Then,   the   conditional   laws   of   $\tau$   given
$\{L_\ca(\tau)=n\}$       and      of       $\tau_{[\theta]}$      given
$\{L_\ca(\tau_{[\theta]})=n\}$ are the same for  all $n\in \N$ such that
$\P\bigl(L_\ca(\tau)=n\bigr) >0$.
\end{prpstn}

Notice that we   don't  assume  that   $p$  is  critical,   sub-critical  or
super-critical in Proposition \ref{prop:prob_equiv}.

\begin{proof}
Let $\bt\in\T_0$. Using \reff{eq:loi-tau} and the definition of 
$p_\theta^\ca$, we have:
\begin{align*}
\P(\tau_{[\theta]}=\bt) &
=\prod_{v\in\bt}p_\theta^\ca\bigl(k_v(\bt)\bigr)\\
&
=\left(\prod_{ k_v(\bt)\in\ca}c_\ca(\theta)
  \theta^{k_v(\bt)-1}p\bigl(k_v(\bt)\bigr)\right)\, 
\left(\prod_{k_v(\bt)\not\in\ca}\theta^{k_v(\bt)-1}
  p\bigl(k_v(\bt)\bigr)\right)\\ 
& 
=c_\ca(\theta)^{L_\ca(\bt)}\theta^{\sum_{v\in\bt}k_v(\bt)
  -\sharp \bt}\, \prod_{v\in\bt}p\bigl(k_v(\bt)\bigr). 
\end{align*}
Since $\sum _{v\in\bt}k_v(\bt)=\sharp \bt-1$, we deduce that:
\begin{equation}
   \label{eq:pq+}
\P(\tau_{[\theta]}=\bt)=\frac{c_\ca(\theta)^{L_\ca(\bt)}}{\theta}\, \P(\tau=\bt).
\end{equation}
By summing \reff{eq:pq+} on $\{\bt\in \T_0, \, L_\ca(\bt )=n\}$ and
using~\eqref{eq:bgw-in-T0}, 
we obtain:
\[
\P\bigl(L_\ca(\tau_{[\theta]})=n\bigr) 
= \frac{c_\ca(\theta)^{n}}{\theta}\P\bigl(L_\ca(\tau)=n\bigr).
\]
Let $n\in \N$ be such that $\P\bigl(L_\ca(\tau)=n\bigr) >0$. 
Since $c_\ca(\theta)$ is positive, by dividing \eqref{eq:pq+} with this equation  term by term, we get that
for $\bt\in \T_0$ such that $L_\ca(\bt)=n$, we have:
\[
\P(\tau=\bt\, |\, L_\ca(\tau)=n)=\P(\tau_{[\theta]}=\bt\, |\, L_\ca(\tau_{[\theta]})=n). 
\]
This ends the proof  as $\tau$ (resp. $\tau_{[\theta]}$) is  a.s.\ finite on 
$\{L_\ca(\tau)=n\}$ (resp. $\{L_\ca(\tau_{[\theta]})=n\}$) by~\eqref{eq:bgw-in-T0}. 
\end{proof}

\begin{rmrk}[Offspring distributions with the  same distribution of
  conditionned BGW]
Let $\tau$ be a BGW tree with offspring distribution $p$  such that:
  \begin{equation}
    \label{eq:degen}
    0<p(0) \quad\text{and}\quad
     p(0)+p(1)<1. 
   \end{equation}
   Let  $\ca\subset  \N$  such  that  $p(\ca)>0$.   It  is  possible  to
   characterize the set $\cp(p, \ca)$ of probability $p'$ satisfying the
   non-degeneracy                            condition~\eqref{eq:degen},
   $\supp(p')\subset  \supp(p)$   and  for  all  $n\in   \N$  such  that
   $\P(L_\ca(\tau')=n)>0$, where $\tau'$ is  the BGW tree with offspring
   distribution   $p'$,   we   have  that   $\tau'$   conditionally   on
   $\{L_\ca(\tau')=n\}$  is  distributed   as  $\tau$  conditionally  on
   $\{L_\ca(\tau)=n\}$.   According to~\cite[Theorem~4.1]{abd:cbgw},  we
   have       that      $\cp(p,       \ca)$       is      equal       to
   $\{p^\ca_\theta\,   \colon\,   \theta\in  \mathring   I_\ca\}$   with
   eventually (depending on  $p$ and $\ca$) one or  two more probability
   distributions associated in some sense  to the parameters $\theta$
   belonging to the    boundary of $I^\ca$ in $[0, +\infty ]$. 
\end{rmrk}
 
\subsubsection{The generic sub-critical case}
\label{sec:gene-sub}

Let $p$ be a sub-critical offspring distribution such that $p(0)+p(1)<1$
and let $\ca\subset\N$ such that  $p(\ca)>0$.  For $\theta\in I_\ca$, we
denote   by  $m^\ca(\theta)$   the   mean  value   of  the   probability
$p_\theta^\ca$. By  hypothesis, we have that  $0<m^\ca(1)<1$.  According
to  Corollary~5.7  and  Proposition~5.10  in~\cite{abd:cbgw},  see  also
\cite{ad:llcgwtcc},  the  continuous  map $\theta  \mapsto  m^\ca(\theta)$
defined on $I^\ca$ is finite on $\mathring I^\ca$, increasing as long as
$m^\ca(\theta)\leq   1$   (and   it   might  be   decreasing   at   some
$\theta\in   I^\ca$    where   $m^\ca(\theta)>1$,   see    Remark   5.12
in~\cite{abd:cbgw}).  We thus deduce the following lemma.

\begin{lmm}[Existence of at most one critical parameter]
\label{lem:generic}
Let   $p$   be   a  sub-critical   offspring   distribution   satisfying
$p(0)+p(1)<1$ and $\ca\subset \N$ such that $p(\ca)>0$. There  exists at
most  one critical value $\theta\in I_\ca$  such that
$m^\ca(\theta)=1$.
\end{lmm}

When it exists, we denote the critical parameter by $\theta_\ca^\mathrm{c}$ given
as the unique solution of
$m^\ca(\theta)=1$ in $I_\ca$;  and we shall consider the critical
offspring distribution: 
\begin{equation}
   \label{eq:def-p*gene}
p^\ca=p_{\theta_\ca^\mathrm{c}}^\ca.
\end{equation}

\begin{dfntn}
The offspring distribution $p$ is said to be generic for the set
$\ca$ if $\theta_\ca^\mathrm{c}$ exists.
\end{dfntn}

Let $d_\ca$ be  the period
of the  offspring distribution associated in
Lemma~\ref{lem:tA} to $p^\ca$ through the Rizzolo's bijection. Notice it
depends only on $\supp(p)$ and $\ca$. 
 Using Proposition \ref{prop:prob_equiv} and Corollary
\ref{cor:kesten-ca}, we immediately deduce the following result in the
sub-critical generic case.

\begin{prpstn}[Generic case: the local limit is a Kesten's tree]
\label{prop:gene-sub}
  Let   $p$  be   a  sub-critical   offspring  distribution  such that
  $p(0)+p(1)<1$  and  let $\ca\subset\N$  such that  $p(\ca)>0$. Assume
  that $p$ is generic for $\ca$.  Let $\tau$ be a BGW tree with offspring
  distribution $p$ and let $\tau^*_\ca$ be a Kesten's tree associated with
  the  offspring distribution  $p^\ca$ given  by~\eqref{eq:def-p*gene}.
 Let $\tau_n$  be a  random tree
  distributed     according      to     $\tau$      conditionally     on
  $\left\{L_\ca(\tau)=1+ d_\ca\,  n\right\}$ (resp.  $\{L_\ca(\tau)\geq
  n\}$).  Then, we have:
\[
\tau_n\; \xrightarrow[n\rightarrow \infty ]{\textbf{(d)}}   \tau^*_\ca. 
\]
\end{prpstn}

\subsubsection{The non-generic sub-critical case}
\label{sec:non-gene}

In order  to state precisely the  general result, we shall  consider the
set $\T_\infty$  of trees that  may have  infinite nodes and  extend the
definition of the local convergence on this set.

Let $n\in \N$. For $u=u_1u_2\ldots u_n\in\cu$, we set
$|u|_\infty=\max\{n,\max\{u_i,1\le i\le n\}\}$ 
and we define the associated restriction operator $r_n^\infty $ on
$T^\infty $ defined by:
\[
\ r_n^\infty (\bt)=\{u\in\bt,\ |u|_\infty\le n\}.
\]
For all tree $\bt\in \T_\infty $, the restricted tree $r_n^\infty (\bt)$
has height at most $n$ and all the nodes have at most $n$ offsprings
(hence the tree $r_n^\infty (\bt)$ is finite).
We define also the associated distance, for all $\bt, \bt'\in \T_\infty
$:
\[
\delta_\infty(\bt,\bt')=2^{-\sup\{n\in\N,\ r_n^\infty (\bt)=r_n^\infty (\bt')\}}.
\]
Remark that, for trees in $\T$,  the topologies induced by the distances
$\delta$ and $\delta_\infty$  coincide. We will from now-on  work on the
space $\T_\infty$ endowed with the distance $\delta_\infty$. It is clear
that the metric space $(\T_\infty,\delta_\infty) $ is compact.

Set $\bar \N=\N \cup\{+\infty \}$. 
If $p=(p(n), n\in \N)$ is a sub-critical offspring distribution with
mean $m<1$, we define   $p^{**}=( p^{**}(n), n\in \bar\N)$ a probability distribution
on $\bar \N$ with:
\[
 p^{**}(n)=np(n)\quad \text{for $n\in \N$ and}\quad 
 p^{**}(+\infty )=1-m.
\] 

We  define a  new  random  tree on  $\T_\infty  $, denoted  $\tau^{**}$,
associated to  $p$ in a way  very similar to the  definition of the Kesten's
tree.

\begin{dfntn}
   \label{defi:c-tree}
Let $p$ be a sub-critical offspring distribution.
The condensation tree $\tau^{**}$ associated with $p$ 
is a two-type BGW tree taking values in $\T_\infty $ and distributed as follows: 
\begin{itemize}
   \item[-] Individuals are \textit{normal} or \textit{special}.
   \item[-] The root of $\tau^{**}$    \textit{special}.
   \item[-] A \textit{normal} individual   produces only
     \textit{normal} individuals
     according to  $p$. 
\item[-] A \textit{special}  individual  produces individuals 
  according to the  distribution $p^{**}$. 
\begin{itemize}
\item[-] If it has a finite number of offsprings,  then one of them 
  chosen uniformly at random, is  
  \textit{special},  the others (if any) are \textit{normal}. 
\item[-] If it has an infinite number of offsprings,  then all of them  are \textit{normal}. 
  \end{itemize}
\end{itemize}
\end{dfntn}

As  we  suppose  that  $p$   is  sub-critical  (that is,   $m<1$),  then  the
condensation  tree $\tau^{**}$  associated  with $p$  has  a.s.\ only  one
infinite node, and its random height  is distributed as $G-1$, where $G$
has the geometric distribution with parameter $1-m$.

Recall definitions from
Section \ref{sec:def-ca}. Set:
\[
  \theta^*_\ca=\sup I_\ca\in [1, +\infty].
\]

The next lemma completes Lemma \ref{lem:generic}, and it is part of Lemma 5.2 in \cite{ad:llcgwtcc}. 
We provide a short proof based on~\cite{abd:cbgw} since $c_\ca(\theta^*_\ca)$ in the definition of $I_\ca$ is assumed to be positive, whereas it is only non-negative in~\cite{ad:llcgwtcc}. 

\begin{lmm}[$\theta^*_\ca$ belongs to $I_\ca$]
\label{lem:non-generic}
Let   $p$   be   a  sub-critical   offspring   distribution   such that
$p(0)+p(1)<1$  and  $\ca\subset \N$  such that  $p(\ca)>0$. If $m^\ca(\theta)<1$ for all
$\theta\in I_\ca$ (that  is, $p$  is not  generic for
$\ca$), then the parameter $\theta^*_\ca$ is finite belongs to $I_\ca$.
\end{lmm}

\begin{proof}
    By continuity, it is possible that $\sum_{k\not\in \ca} (\theta^*_\ca)^{k-1} p(k)=1$ and 
$\sum_{k\in \ca} (\theta^*_\ca)^{k-1} p(k)<+\infty$ (and thus $\theta^*_\ca$ belongs to $(1,+\infty)$). This corresponds to the case $\beta=0$ in \cite[Remark 3.7.d]{abd:cbgw} (with $J=1$) and thus $m^\ca(\theta^*_\ca)>1$. So assuming that $m^\ca<1$ in  $I_\ca$ implies by monotone convergence that $\sum_{k\in \ca} (\theta^*_\ca)^{k-1} p(k)$ is finite and 
thus that $\sum_{k\not\in \ca} (\theta^*_\ca)^{k-1} p(k)<1$, that is, $\theta^*_\ca$ belongs to $I_\ca$. 
\end{proof}

When $p$  is not  generic for
$\ca$, we shall consider the sub-critical offspring distribution:
\[
p^{\ca}=p_{\theta_\ca^*}^\ca.
\]
Let $d_\ca$ denote the period
of the  offspring distribution  associated in
Lemma~\ref{lem:tA} to $p^{\ca}$ through Rizzolo's bijection. Notice it
depends only on $\supp(p)$ and $\ca$ since $\supp(p^\ca)=\supp(p)$. 

Following the idea  developed for the critical case,  with more involved
technicalities, we can prove the following result, see
\cite[Theorem~1.3]{ad:llcgwtcc}; this extends results on condensation from~\cite{j:sgtcgwrac,js:cnt} 
where only the case $\ca=\N$ (i.e. the total population size) was considered.

\begin{prpstn}[Non-generic case:  the local limit is a
  condensation  tree]
\label{prop:non-gene-sub}
Let $p$ be a sub-critical offspring distribution such that $p(0)+p(1)<1$  and let
$\ca\subset\N$ such that $p(\ca)>0$. Assume that  $p$ is not generic
for $\ca$.
Let $\tau$ be a BGW
tree with offspring distribution $p$ and let $ \tau^{**}_\ca$ be a
condensation  tree associated with
the sub-critical offspring distribution $ p ^ \ca$. 
Let $\tau_n$ be a random
  tree    distributed    according    to   $\tau$    conditionally    on
  $\left\{L_\ca(\tau)=1+ d_A\, n\right\}$ (resp.  $\{L_\ca(\tau)\geq n\}$). 
Then, we have: 
\[
\tau_n\; \xrightarrow[n\rightarrow \infty ]{\textbf{(d)}}    \tau^{**}_\ca. 
\]
\end{prpstn}

\begin{rmrk}
Considering the local limit of BGW aims at understanding the shape of the random tree near the root. 
The local limit of the BGW tree rerooted at a random node is studied in~\cite{stufler19} with $\ca=\N$ in the generic and non-generic cases; in both cases the limiting tree is similar either to the Kesten's tree or the condensation tree. 

The condensation tree appears also as the limit of subcritical BGW with different conditioning: see~\cite{h:cgwtmod} where the BGW tree is conditioned on the maximal out-degree to be large or~\cite{sonia} where the nodes of the tree are marked independently with a probability depending on their out-degree and the tree is conditioned to have a large number of marks.
\end{rmrk}

\subsubsection{Generic and non-generic distributions}
\label{sec:gene-non-gene}
Let   $p$   be   a  sub-critical   offspring   distribution such that
$p(0)+p(1)<1$.  We shall give a criterion to say   easily for which sets $\ca$ the
offspring distribution $ p$ is  generic. As we have $m<1$, we  want to find a
$\theta$    (which    will    be    greater   than    1)    such    that
$m^\ca(\theta)=1$. This  problem is closely related  to the radius $\rho\geq
1$ of convergence of the generating function of $p$, denoted by $g$.

We have the following result.

\begin{lmm}[\cite{ad:llcgwtcc}, Lemma 5.4]
\label{lem:non-gene}
Let $p$ be a sub-critical offspring distribution such that
$p(0)+p(1)<1$.
\begin{itemize}
\item[(i)] If $\rho=+\infty$ or if ($\rho<+\infty$ and $g'(\rho)\ge 1$),
  then $p$ is generic for any  $\ca\subset \N$ such that $p(\ca)>0$.
\item[(ii)]  If $\rho=1$ (that is,  the probability $p$ admits no exponential
  moment), then $p$ is non-generic for every $\ca\subset \N$ such that $p(\ca)>0$.
\item[(iii)]  If $1<\rho <+\infty$ and $g'(\rho)<1$, then $p$ is
  non-generic for $\ca\subset \N$, with $p(\ca)>0$,  if and only if:
\[
\E[Y|Y\in\ca]<\frac{\rho-\rho g'(\rho)}{\rho-g(\rho)},
\]
where  $Y$ is distributed according to $p_\rho^\N$, that is:
\[
  \E[f(Y)]=\frac{\E[f(\zeta) \rho^\zeta]}{g(\rho)}\cdot
\] 
In particular,  $p$ is non-generic for $\ca=\{0\}$ but
generic of $\ca=\{k\}$ for any $k$ large enough such that $p(k)>0$.
\end{itemize}
\end{lmm}

\begin{rmrk}
\label{rem:gene}
  In  case (iii)  of  Lemma \ref{lem:non-gene},  we  gave in Remark 5.5 of 
\cite{ad:llcgwtcc}:
\begin{itemize}
   \item a sub-critical offspring
  distribution which is generic for  $\N$ but non-generic for $\{0\}$;
\item a sub-critical offspring
  distribution which is non-generic for  $\N$ but
  generic  for $\{k\}$  with  $k$  large enough. 
\end{itemize}
This   shows that  the
  genericity of sets is not monotone with respect to the inclusion.
\end{rmrk}

\subsection{Other results}\label{sec:extension}
We mention some related problem. 
\begin{enumerate}
  
\item \textbf{Marked BGW tree}. 
The proof of Theorem \ref{theo:main1} can be slighly modified to study
a tree with randomly marked nodes: conditionally given the tree, we
mark its nodes randomly, independently of each others, with a
probability that depends only on the out-degree of the node. Then we
obtain that a critical BGW tree conditioned on having $n$ marked nodes
still converges in distribution toward a Kesten's tree, see
\cite{abd:llgwtcnpn}. This allows to study a Galton-Watson tree
conditioned on the number of protected nodes where a protected node is
a node that is neither a leaf nor the parent of a leaf.
See also \cite{sonia} for the sub-critical case in the generic and non-generic case.

\item \textbf{Very fat BGW trees}. 
  Recall  that  $Z_n=z_n  (\tau)$  represents  the  population  size  at
  generation $n$  of the BGW  tree $\tau$.  In  the same spirit as in Section~\ref{sec:large-gen},  we can
  consider a critical BGW tree  $\tau$ conditioned on $\{Z_n=a_n\}$ with
  a sequence $(a_n, n\in \N^*)$  of positive integers.  Notice this does
  not  fall  into the  framework  developed  in  Section~\ref{sec:criteria} for  a
  functional   $A$   with   \ref{eq:>}-\ref{eq:+}-\ref{eq:=}   property.
  Nevertheless, the  proof of  Theorem \ref{theo:main1} can  be adapted,
  see~\cite{ad:llcgwtisc},  to  prove  that  a critical  BGW  tree  with
  geometric  offspring distribution  conditioned  on $\{Z_n=a_n\}$  with
  $\lim_{n\rightarrow  \infty  }  a_n/n^2=0$ converges  in  distribution
  toward the Kesten's tree.  The case  $a_n\sim C n^2$ is more involved as
  the infinite spine is replaced by an infinite backbone that  does not satisfy the usual  branching property: the
  numbers  of offspring  inside a  generation are  not independent,  see
  \cite{abd:nlcgwt} for  details.  This new  limit appears  also when  considering large  BGW
   trees with $n$ nodes and exponential  weight given by its height, see~\cite{dm:exp-w}. 
   
   Let  us mention  that there  are some
  results, with the same  flavor, established in~\cite{ad:apegw} for the
  super-critical  case  that  can  also  be adapted  to  BWG  tree  with
  sub-critical offspring distribution $p$ such  that $g(r)=r$ has a root
  strictly  larger  than 1  (this  amounts  for $p$  to  be  of the  form
  $\tilde  p$  in~\eqref{eq:tilde-p}).
  The  local  convergence  of  the  critical  BGW  tree
  conditioned to $\{Z_n=a_n\}$  with a sequence $(a_n,  n\in \N^*)$ such
  that $\lim_{n\rightarrow \infty } a_n=+\infty $ is an open question in
  general.

  \item\textbf{Multi-type BGW tree}.
The ideas from Section~\ref{sec:out-d} for the conditioning on $\{L_\ca(\tau)=n\}$ in the critical case
can also be studied for  critical 
multi-type BGW trees conditioned on the size the population of each type, see \cite{adg:cmtgwtcl} where the limit is now 
a multi-type Kesten's tree. In this setting, the asymptotic proportion of each type is fixed equal to the normalized left-eigen-vector of the mean matrix associated to the eigen-value 1.
See also \cite{p:vwcmgwp,s:lclcmgwtarm,th:large-GW} for
other results on this topic. 
Let us mention that conditioning a BGW tree $\tau$ on the event $\{L_{\ca_1}(\tau)(\tau)=n_1, \ldots, L_{\ca_k}(\tau)=n_k\}$, where the sets $\ca_1, \ldots, \ca_k$ are pairwise disjoint, and $n_1+ \ldots +n_k$ goes to infinity, is studied in~\cite{abd:cbgw} in  the generic case using a multi-type BGW tree approach. The non generic case for this problem and more generally for conditioned multi-type sub-critical BGW tree is still an open question.

\item \textbf{L\'evy tree}. 
   There is  also a natural extension of this work when considering
scaling limit of conditioned critical BGW:
   see  \cite{d:ltcpcgwt, kor:GW-cond-leaf} when conditioning on the total size or the number of leaves, see also
\cite{kortchemski2015condensation} for scaling limits in the sub-critical non-generic case (where condensation occurs). 
It is also possible to consider the same conditioning directly on the L\'evy tree which appears as the scaling limit of the BGW tree:  see~\cite{xin-loc-cv,adn:jump} when conditioning on a node to be
   very large or~\cite{adh:large,n:root} when conditioning on the total size to be large. 

\end{enumerate}

\bibliographystyle{abbrv}
\bibliography{biblio}

\begin{thebibliography}{10}

\bibitem{abd:cbgw}
R.~ABRAHAM, H.~BI, and J.-F. DELMAS.
\newblock Conditioning {B}ienaym\'e-{G}alton-{W}atson trees to have large
  sub-populations.
\newblock arXiv:2311.17716, 2023.

\bibitem{abd:llgwtcnpn}
R.~ABRAHAM, A.~BOUAZIZ, and J.-F. DELMAS.
\newblock Local limits of {G}alton-{W}atson trees conditioned on the number of
  protected nodes.
\newblock {\em J. of Appl. Probab.}, 54:55--65, 2017.

\bibitem{abd:nlcgwt}
R.~ABRAHAM, A.~BOUAZIZ, and J.-F. DELMAS.
\newblock Very fat geometric {G}alton-{W}atson trees.
\newblock {\em ESAIM: PS}, 24:294--314, 2020.

\bibitem{sonia}
R.~ABRAHAM, S.~BOULAL, and P.~DEBS.
\newblock Condensation for sub-crtical marked {G}alton-{W}atson tree.
\newblock Work in progress.

\bibitem{ad:llcgwtcc}
R.~ABRAHAM and J.-F. DELMAS.
\newblock Local limits of conditioned {G}alton-{W}atson trees: the condensation
  case.
\newblock {\em Elec. J. of Probab.}, 19:Article 56, 1--29, 2014.

\bibitem{ad:llcgwtisc}
R.~ABRAHAM and J.-F. DELMAS.
\newblock Local limits of conditioned {G}alton-{W}atson trees: the infinite
  spine case.
\newblock {\em Elec. J. of Probab.}, 19:Article 2, 1--19, 2014.

\bibitem{ad:apegw}
R.~ABRAHAM and J.-F. DELMAS.
\newblock Asymptotic properties of expansive {G}alton-{W}atson trees.
\newblock {\em Electron. J. Probab.}, 24:Paper No. 15, 51, 2019.

\bibitem{adg:cmtgwtcl}
R.~ABRAHAM, J.-F. DELMAS, and H.~GUO.
\newblock Critical multi-type galton-watson trees conditioned to be large.
\newblock {\em J. of Th. Probab.}, 31:757--788, 2018.

\bibitem{adh:pgwttvmp}
R.~ABRAHAM, J.-F. DELMAS, and H.~HE.
\newblock Pruning {G}alton-{W}atson trees and tree-valued {M}arkov processes.
\newblock {\em Ann. de l'Inst. Henri Poincar\'e}, 48:688--705, 2012.

\bibitem{adh:large}
R.~ABRAHAM, J.-F. DELMAS, and H.~HE.
\newblock Brownian continuum random tree conditioned to be large.
\newblock arXiv:2202.10258, 2022.

\bibitem{adn:jump}
R.~ABRAHAM, J.-F. DELMAS, and M.~NASSIF.
\newblock Conditioning (sub)critical l{\'e}vy trees by their maximal degree:
  Decomposition and local limit.
\newblock arXiv:2211.02317, 2022.

\bibitem{a:crt3}
D.~ALDOUS.
\newblock The continuum random tree {III}.
\newblock {\em Ann. of Probab.}, 21:248--289, 1993.

\bibitem{ap:tvmcdgwp}
D.~ALDOUS and J.~PITMAN.
\newblock Tree-valued {M}arkov chains derived from {G}alton-{W}atson processes.
\newblock {\em Ann. de l'Inst. Henri Poincar\'e}, 34:637--686, 1998.

\bibitem{an:bp}
K.~ATHREYA and P.~NEY.
\newblock {\em Branching processes}.
\newblock Springer-Verlag, New York-Heidelberg, 1972.
\newblock Die Grundlehren der mathematischen Wissenschaften, Band 196.

\bibitem{b:gw}
I.~BIENAYME.
\newblock De la loi de multiplication et de la dur{\'e}e des familles.
\newblock {\em Soc. Philomath. Paris Extraits}, S\'er. 5:37--39, 1845.

\bibitem{b:cpm}
P.~BILLINGSLEY.
\newblock {\em Convergence of probability measures}.
\newblock Wiley Series in Probability and Statistics: Probability and
  Statistics. John Wiley \& Sons Inc., New York, second edition, 1999.
\newblock A Wiley-Interscience Publication.

\bibitem{bdr:HS21}
A.~BRANDENBERGER, L.~DEVROYE, and T.~REDDAD.
\newblock {The Horton-Strahler number of conditioned {G}alton-{W}atson trees}.
\newblock {\em Electronic Journal of Probability}, 26:1 -- 29, 2021.

\bibitem{bru}
B.~BRU.
\newblock {\`A{}} la recherche de la d\'emonstration perdue de {B}ienaym\'e.
\newblock {\em Math. Inform. Sci. Humaines}, 114:5--17, 1991.

\bibitem{ck:rncpccgwta}
N.~CURIEN and I.~KORTCHEMSKI.
\newblock Random non-crossing plane configurations: a conditioned
  {G}alton-{W}atson tree approach.
\newblock {\em Random Struct. and Alg.}, 45:236--260, 2014.

\bibitem{d:rt}
M.~DRMOTA.
\newblock {\em Random trees}.
\newblock SpringerWienNewYork, Vienna, 2009.
\newblock An interplay between combinatorics and probability.

\bibitem{dp:rftt}
M.~DRMOTA and H.~PRODINGER.
\newblock The register function for {$t$}-ary trees.
\newblock {\em ACM Trans. Algorithms}, 2(3):318--334, 2006.

\bibitem{d:ltcpcgwt}
T.~DUQUESNE.
\newblock A limit theorem for the contour process of conditioned
  {G}alton-{W}atson trees.
\newblock {\em Ann. of Probab.}, 31:996--1027, 2003.

\bibitem{dlg:rt}
T.~DUQUESNE and J.-F. LE~GALL.
\newblock {\em Random trees, {L}\'evy processes and spatial branching
  processes}.
\newblock Number 281. Ast\'erisque, 2002.

\bibitem{dm:exp-w}
B.~DURHUUS and M.~\"UNEL.
\newblock Trees with exponential height dependent weight.
\newblock {\em Probab. Theory Related Fields}, 186(3-4):999--1043, 2023.

\bibitem{d:pte}
R.~DURRETT.
\newblock {\em Probability: theory and examples}.
\newblock Duxbury Press, Belmont, CA, second edition, 1996.

\bibitem{d:tpbprrw}
M.~DWASS.
\newblock The total progeny in a branching process and a related random walk.
\newblock {\em J. Appl. Probability}, 6, 1969.

\bibitem{et:cv-tree}
G.~ELEK and G.~TARDOS.
\newblock Convergence and limits of finite trees.
\newblock {\em Combinatorica}, 42(6):821--852, 2022.

\bibitem{esparza}
J.~ESPARZA, M.~LUTTENBERGER, and M.~SCHLUND.
\newblock A brief history of {S}trahler numbers.
\newblock In {\em Language and automata theory and applications}, volume 8370
  of {\em Lecture Notes in Comput. Sci.}, pages 1--13. Springer, Cham, 2014.

\bibitem{e:rt}
S.~N. EVANS.
\newblock {\em Probability and real trees}, volume 1920 of {\em Lecture Notes
  in Mathematics}.
\newblock Springer, Berlin, 2008.
\newblock Lectures from the 35th Summer School on Probability Theory held in
  Saint-Flour, July 6--23, 2005.

\bibitem{gw:pef}
F.~GALTON and H.~WATSON.
\newblock On the probability of extinction of families.
\newblock {\em J. Anthropol. Inst.}, 4:138--144, 1874.

\bibitem{haas2012}
B.~HAAS and G.~MIERMONT.
\newblock Scaling limits of {M}arkov branching trees with applications to
  {G}alton-{W}atson and random unordered trees.
\newblock {\em Ann. Probab.}, 40(6):2589--2666, 2012.

\bibitem{hjv:pb}
P.~HACCOU, P.~JAGERS, and V.~A. VATUTIN.
\newblock {\em Branching processes: variation, growth, and extinction of
  populations}, volume~5 of {\em Cambridge Studies in Adaptive Dynamics}.
\newblock Cambridge University Press, Cambridge; IIASA, Laxenburg, 2007.

\bibitem{h:tbp}
T.~E. HARRIS.
\newblock {\em The theory of branching processes}.
\newblock Die Grundlehren der Mathematischen Wissenschaften, Bd. 119.
  Springer-Verlag, Berlin; Prentice-Hall, Inc., Englewood Cliffs, N.J., 1963.

\bibitem{h:cgwtmod}
X.~HE.
\newblock Conditioning {G}alton-{W}atson trees on large maximal out-degree.
\newblock {\em J. of Th. Probab.}, 30:842--851, 2017.

\bibitem{xin-loc-cv}
X.~HE.
\newblock Local convergence of critical random trees and continuous-state
  branching processes.
\newblock {\em J. Theoret. Probab.}, 35(2):685--713, 2022.

\bibitem{j:sgtcgwrac}
S.~JANSON.
\newblock Simply generated trees, conditioned {G}alton-{W}atson trees, random
  allocations and condensation.
\newblock {\em Probab. Surv.}, 9:103--252, 2012.

\bibitem{janson:fringe}
S.~JANSON.
\newblock Asymptotic normality of fringe subtrees and additive functionals in
  conditioned {G}alton-{W}atson trees.
\newblock {\em Random Structures Algorithms}, 48(1):57--101, 2016.

\bibitem{js:cnt}
T.~JONNSSON and S.~STEFANSSON.
\newblock Condensation in nongeneric trees.
\newblock {\em J. Stat. Phys.}, 142:277--313, 2011.

\bibitem{k:history}
D.~G. KENDALL.
\newblock The genealogy of genealogy: branching processes before (and after)
  1873.
\newblock {\em Bull. London Math. Soc.}, 7(3):225--253, 1975.

\bibitem{k:gwctp}
D.~KENNEDY.
\newblock The {G}alton-{W}atson process conditioned on the total progeny.
\newblock {\em J. Appl. Probability}, 12(4):800--806, 1975.

\bibitem{k:sbrwrc}
H.~KESTEN.
\newblock Subdiffusive behavior of random walk on a random cluster.
\newblock {\em Ann. Inst. H. Poincar\'e Probab. Statist.}, 22(4):425--487,
  1986.

\bibitem{ks:ltmgwp}
H.~KESTEN and B.~P. STIGUM.
\newblock A limit theorem for multidimensional {G}alton-{W}atson processes.
\newblock {\em Ann. Math. Statist.}, 37:1211--1223, 1966.

\bibitem{khanfir2024}
R.~KHANFIR.
\newblock Fluctuations of the {H}orton-{S}trahler number of stable
  {G}alton-{W}atson trees.
\newblock arXiv:2401.13771, 2024.

\bibitem{ka:bpb}
M.~KIMMEL and D.~E. AXELROD.
\newblock {\em Branching processes in biology}, volume~19 of {\em
  Interdisciplinary Applied Mathematics}.
\newblock Springer-Verlag, New York, 2002.

\bibitem{kor:GW-cond-leaf}
I.~KORTCHEMSKI.
\newblock Invariance principles for {G}alton-{W}atson trees conditioned on the
  number of leaves.
\newblock {\em Stochastic Process. Appl.}, 122(9):3126--3172, 2012.

\bibitem{kortchemski2015condensation}
I.~KORTCHEMSKI.
\newblock Limit theorems for conditioned non-generic {G}alton-{W}atson trees.
\newblock {\em Ann. Inst. Henri Poincar\'{e} Probab. Stat.}, 51(2):489--511,
  2015.

\bibitem{kr:cond-cauchy}
I.~KORTCHEMSKI and L.~RICHIER.
\newblock Condensation in critical {C}auchy {B}ienaym\'e-{G}alton-{W}atson
  trees.
\newblock {\em Ann. Appl. Probab.}, 29(3):1837--1877, 2019.

\bibitem{lpp:cp}
R.~LYONS, R.~PEMANTLE, and Y.~PERES.
\newblock Conceptual proofs of {$L\log L$} criteria for mean behavior of
  branching processes.
\newblock {\em Ann. Probab.}, 23(3):1125--1138, 1995.

\bibitem{m:nvgdgwt}
N.~MIMAMI.
\newblock On the number of vertices with a given degree in a {G}alton-{W}atson
  tree.
\newblock {\em Adv. in Appl. Probab.}, 37(1):229--264, 2005.

\bibitem{n:root}
M.~NASSIF.
\newblock Zooming in at the root of the stable tree.
\newblock {\em Electron. J. Probab.}, 27:Paper No. 39, 38, 2022.

\bibitem{n:tece}
J.~NEVEU.
\newblock Sur le th\'eor\`eme ergodique de {C}hung-{E}rd{\H o}s.
\newblock {\em C. R. Acad. Sci. Paris}, 257:2953--2955, 1963.

\bibitem{n:apghw}
J.~NEVEU.
\newblock Arbres et processus de {G}alton-{W}atson.
\newblock {\em Ann. de l'Inst. Henri Poincar\'e}, 22:199--207, 1986.

\bibitem{p:vwcmgwp}
S.~PENISSON.
\newblock Beyond the q-process: various ways of conditioning multitype
  {G}alton-{W}atson processes.
\newblock {\em ALEA, Lat. Am. J. Probab. Math. Stat.}, 13:223--237, 2016.

\bibitem{p:csp}
J.~PITMAN.
\newblock {\em Combinatorial stochastic processes}, volume 1875 of {\em Lecture
  Notes in Mathematics}.
\newblock Springer-Verlag, Berlin, 2006.
\newblock Lectures from the 32nd Summer School on Probability Theory held in
  Saint-Flour, July 7--24, 2002.

\bibitem{r:slmbtgwtcnvodgs}
D.~RIZZOLO.
\newblock Scaling limits of {M}arkov branching trees and {G}alton-{W}atson
  trees conditioned on the number of vertices with out-degree in a given set.
\newblock {\em Ann. de l'Inst. Henri Poincar\'e}, 51(2):512--532, 2015.

\bibitem{s:prw}
F.~SPITZER.
\newblock {\em Principles of random walk}.
\newblock Springer-Verlag, New York, second edition, 1976.
\newblock Graduate Texts in Mathematics, Vol. 34.

\bibitem{s:lclcmgwtarm}
R.~STEPHENSON.
\newblock Local convergence of large critical multi-type {G}alton-{W}atson
  trees and applications to random maps.
\newblock {\em J. of Th. Probab.}, 31:159--205, 2018.

\bibitem{stufler19}
B.~STUFLER.
\newblock Local limits of large {G}alton-{W}atson trees rerooted at a random
  vertex.
\newblock {\em Ann. Inst. Henri Poincar\'e{} Probab. Stat.}, 55(1):155--183,
  2019.

\bibitem{stufler20}
B.~STUFLER.
\newblock {On the maximal offspring in a subcritical branching process}.
\newblock {\em Electronic Journal of Probability}, 25:1 -- 62, 2020.

\bibitem{th:large-GW}
P.~TH\'EVENIN.
\newblock Vertices with fixed outdegrees in large {G}alton-{W}atson trees.
\newblock {\em Electron. J. Probab.}, 25:Paper No. 64, 25, 2020.

\bibitem{viennot}
G.~X. VIENNOT.
\newblock Trees.
\newblock In {\em Mots}, Lang. Raison. Calc., pages 265--297. Herm\`es, Paris,
  1990.

\end{thebibliography}

\end{document}